\documentclass[a4paper,11pt]{article}
\usepackage{amsmath, amsthm, amssymb}
\usepackage{picture}

\usepackage{color}

\newtheorem{thm}{Theorem}[section]
\newtheorem{pro}[thm]{Proposition}
\newtheorem{lem}[thm]{Lemma}

\theoremstyle{definition}
\newtheorem{defn}[thm]{Definition}
\theoremstyle{remark}
\newtheorem{rem}[thm]{Remark}

\newtheorem{exa}[thm]{Example}

\newenvironment{prf}{\noindent \textbf{Proof:} \ }{\hfill $\Box$}

\newcommand{\CC}{\mathbb{C}}

\newcommand{\PP}{\mathbb{P}}
\newcommand{\QQ}{\mathbb{Q}}
\newcommand{\ZZ}{\mathbb{Z}}
\newcommand{\RR}{\mathbb{R}}

\newcommand{\Ocal}{\mathcal{O}}

\newcommand{\Lcal}{\mathcal{L}}
\newcommand{\Hcal}{\mathcal{H}}
\newcommand{\Fcal}{\mathcal{F}}
\newcommand{\Wcal}{\mathcal{W}}
\newcommand{\Mcal}{\mathcal{M}}
\newcommand{\Dcal}{\mathcal{D}}

\newcommand{\Aut}{\operatorname{Aut}}
\newcommand{\Ker}{\operatorname{Ker}}

\newcommand{\ad}{\operatorname{ad}}
\newcommand{\Res}{\operatorname{Res}}
\newcommand{\res}{\operatorname{res}}
\newcommand{\can}{\operatorname{can}}
\newcommand{\Id}{\operatorname{Id}}

\newcommand{\al}{\alpha}

\newcommand{\p}{\partial}
\newcommand{\g}{\mathfrak{g}}
\newcommand{\h}{\mathfrak{h}}
\newcommand{\fb}{\mathfrak{b}}
\newcommand{\fn}{\mathfrak{n}}
\newcommand{\s}{\mathfrak{s}}
\newcommand{\hg}{\hat{\g}}

\begin{document}

\title{
\textbf{BCFG Drinfeld-Sokolov Hierarchies and FJRW-Theory}}
\author{Si-Qi Liu, Yongbin Ruan and Youjin Zhang}
\maketitle

\tableofcontents

\section{Introduction}\label{introd}

In 1991, Witten \cite{Wi1}  proposed a remarkable conjecture relating the intersection theory of the Deligne-Mumford moduli space 
$\overline{\Mcal}_{g,k}$ to the Kortweg-de Vries (KdV) hierarchy. 
The
geometric side of the Witten conjecture concerns a complicated generating 
function (or free energy)
$$\Fcal(\hbar, t_0,t_1,\cdots)=\sum_{g\geq 0}\hbar^{g-1}\sum_{k\geq 0}\frac{t_{l_1}\cdots
t_{l_k}}{k!}\langle \tau_{l_1}, \tau_{l_2}, \cdots,
\tau_{l_k}\rangle_g$$
of certain intersection numbers $\langle \tau_{l_1}, \tau_{l_2}, \cdots,
\tau_{l_k}\rangle_g$ on $\overline{\Mcal}_{g,k}$. 
The integrable hierarchy side of the Witten conjecture is a hierarchy of evolutionary PDEs of the form
 $$\frac{\partial u}{\partial t_n}=R_n(u, u_x, u_{xx}, \cdots)$$
for a function of infinitely many time variables $u(x, t_1, t_2, \cdots)$, where $x$ is
a spatial variable and $t_1, t_2, \cdots,  $ are time variables, and $R_n$ are certain polynomials.  The Witten conjecture 
asserts that the generating function $\Fcal$ of the intersection numbers on the Deligne-Mumford moduli space yields a solution of the KdV hierarchy. More precisely, the  total descendant potential function
(or  the partition function) $\Dcal=e^{\Fcal}$
is a tau function of the KdV hierarchy. 

The Witten conjecture was soon proved by
Kontsevich \cite{Ko}. Since then, the Witten-Kontsevich theorem
has brought together the two seemingly alien subjects of integrable
hierarchies and geometry. The intersection theory of the Deligne-Mumford moduli space can be viewed as the Gromov-Witten
theory of a point. 
Immediately after, a great deal of
effort was spent in investigating other integrable hierarchies in
Gromov-Witten theory. A much-studied example is the 2-Toda hierarchy
for the equivariant Gromov-Witten theory of $\PP^1$ by Okounkov-Pandharipande \cite{OP}. Moreover, the Gromov-Witten theory of $\PP^1$ satisfies  a certain reduction of the 2-Toda hierarchy, the so-called extended Toda hierarchy
 by Dubrovin-Zhang and Getzler \cite{DZ4, Ge, Z}. It was generalized later
to orbifold $\PP^1$ \cite{J, MT, PR}. Gromov-Witten theory is a part of the so-called A-model. In many ways, the integrable
hierarchies can be treated as a part of the B-model. We refer to this type of problem as
{\bf Integrable Hierarchy Mirror Symmetries}. 
 The common
characteristics of integrable hierarchy mirror symmetries are: (1) all of them are very
difficult; (2) all of them are mysterious. In particular, the
choice of integrable hierarchies seems to be a matter of luck and
there is no general pattern to predict the integrable hierarchy for a given
A-model geometry.
 
 Instead, one can start from an integrable hierarchy and try to match the A-model geometry with the given integrable
 hierarchies. From this point of view, a natural class of integrable hierarchies is the Drinfeld-Sokolov hierarchies \cite{DS}, of which the KdV-hierachy is the simplest example. A fundamental problem is to identify the A-model geometries which the Drinfeld-Sokolov
 hierarchies govern. This question was very much in Witten's mind when he
proposed his famous conjecture in the first place. Recall that the Drinfeld-Sokolov hierarchies are indexed by
affine Kac-Moody algebras. In particular, we call the integrable hierarchies that correspond to the simply laced affine Kac-Moody algebras the ADE Drinfeld-Sokolov hierarchies. Around the same
time as his conjecture for the KdV hierarchy, Witten proposed a sweeping generalization of his conjecture
\cite{Wi2, Wi3} to identify the A-model geometry of the ADE Drinfeld-Sokolov hierarchies.  The core of his generalization is a remarkable
first order nonlinear elliptic PDE associated to an arbitrary
quasihomogeneous singularity. During the last decade, Witten's
generalization has been explored and a new Gromov-Witten type
theory (FJRW-theory) has been constructed by Fan-Jarvis-Ruan \cite{FJR1,
FJR2}. In particular, Witten's conjecture for the ADE integrable
hierarchies has been verified for the $A_n$ case by Faber-Shadrin-Zvonkine \cite{FSZ} and the $DE$ case by Fan-Jarvis-Ruan \cite{FJR2}. 
We should mention that Witten's original conjecture needs a correction to account for the appearance of mirror symmetry phenomena. This is captured by the following theorem of Fan-Jarvis-Ruan (Theorem 6.1.3 of \cite{FJR2}):
\begin{thm}
\begin{description}
\item[(A)] The total descendant potential function ${\cal D}_{W, G_{max}}$ for the FJRW-theory of $D^T_n, A_{n}, E_6, E_7, E_8$ polynomials $W$ with the
maximal diagonal symmetry group $G_{max}$ is a tau function of the $W^T$ Drinfeld-Sokolov hierarchy.
\item[(B)] For $D_{2n}$, the symmetry group generated by $J=exp(2\pi i \frac{1}{2n-1}, 2\pi i \frac{n-1}{2n-1})$ is different from the maximal diagonal symmetry group. 
The total descendant potential function ${\cal D}_{D_{2n}, \langle J\rangle}$ of FJRW-theory is a tau function of the $D_{2n}$ Drinfeld-Sokolov hierarchy.
\end{description}
\end{thm}

The polynomials $W$ of the above theorem are given by
\[
\begin{array}{lll}
A_n: \  x^{n+1}+y^2, &\quad D_n: \  x^{n-1}+xy^2,&\quad D^T_n: \  x^{n-1}y+y^2, \\
E_6: \  x^3+y^4, &\quad E_7: \  x^3+xy^3,  &\quad E_8: \  x^3+y^5.
\end{array}
\]
For a so-called invertible quasi-homogeneous polynomial $W$, $W^T$ is the ``mirror", or ``transpose", singularity. $A_n$ and
$E_{6,7,8}$ are self-mirror while $D_n$ is mirror to $D^T_n$, as the notation suggests.
For a general invertible quasi-homogeneous polynomial the mirror symmetry  phenomena was studied in \cite{Kr}. 

There are many more Drinfeld-Sokolov hierarchies beyond the ADE cases. The next natural classes are the
$B_n, C_n, F_4, G_2$ Drinfeld-Sokolov 
hierarchies. Unfortunately, the situation becomes much more subtle. To explain the new subtlety, we need to go back to Fan-Jarvis-Ruan's proof of the above
theorem. Their basic idea is that FJRW-theory should be considered as the so-called Landau-Ginzburg A-model. The Landau-Ginzburg B-model is the miniversal deformation  of ADE singularity, which 
carries a formal genus zero Gromov-Witten type theory known as a Frobenius manifold by Saito \cite{S2} and a higher genus theory by Givental \cite{Giv1, Giv2} and Dubrovin-Zhang \cite{DZ3}.  Here, the notion of Frobenius manifold was introduced by Dubrovin  in \cite{Du1, Du2}. 
Note that Saito's theory is  semisimple (see the definition in Section \ref{sec-21}). Dubrovin-Zhang \cite{DZ3} have built an axiomatic theory of integrable
hierarchies for semisimple Frobenius manifolds. In particular, they proved that the higher genus free energies are uniquely determined by the genus zero free energy and the 
so-called Virasoro constraints (see the Appendix), and that the partition function gives a tau function of 
the associated integrable hierarchy. 
This tau function coincides with the total descendant potential defined for a semisimple Frobenius manifold by Givental.
We shall  take the Saito-Givental-Dubrovin-Zhang theory as the definition of the Landau-Ginzburg B-model.
 We should mention that the Givental-Dubrovin-Zhang higher genus theory is only defined over semisimple loci. The extension to non-semisimple loci is a well-known and difficult problem, which has been resolved recently by Milanov \cite{Mi} in the LG-setting. 
 
 Fan-Jarvis-Ruan's proof has two steps. The first step is a higher genus
mirror symmetry theorem to identify the FJRW-theory of $(W, G_{max})$ with the Saito-Givental-Dubrovin-Zhang theory of $W^T$. An early theorem of Frenkel-Givental-Milanov \cite{GM, FGM, Wu1} showed
that Saito-Givental-Dubrovin-Zhang's total descendant potential function is a tau function of the so called Kac-Wakimato hierarchy \cite{KW}. A rather non-trivial theorem of \cite{HM, Wu2} connected Drinfeld-Sokolov hierarchies with Kac-Wakimoto hierarchies. 

There is a Saito-Givental-Dubrovin-Zhang theory for the $B_n$, $C_n$, $F_4$, $G_2$ singularities as well. A natural
expectation is that their total descendant potential functions should be tau functions of 
the corresponding Drinfeld-Sokolov hierarchies. Unfortunately, the following
theorem of Dubrovin-Liu-Zhang \cite{DLZ2} (see Subsection \ref{jw-12-2}) gave a negative answer.

\begin{thm}\label{jw-3} 
The Saito-Givental-Dubrovin-Zhang's total descendant potential functions  of $B_n, C_n, F_4, G_2$ singularities  are NOT tau functions of the corresponding Drinfeld-Sokolov hierarchies.
\end{thm}

Therefore, a new idea is needed to identify the geometry of  $B_n, C_n, F_4, G_2$ Drinfeld-Sokolov hierarchies. This is the main goal of this paper. 

Recall that the $B_n, C_n, F_4, G_2$ series
can be obtained by folding the following Dynkin diagrams of ADE singularities using their
finite symmetry groups:

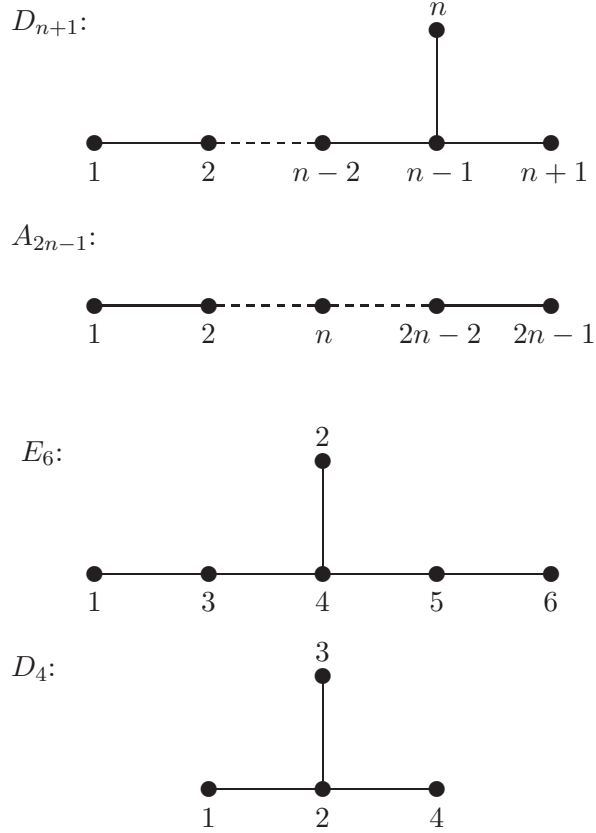
\begin{center}
\setlength{\unitlength}{1mm}
\begin{picture}(86,24)
\put(11,5){\circle*{2}}
\put(26,5){\circle*{2}}
\put(41,5){\circle*{2}}
\put(56,5){\circle*{2}}
\put(71,5){\circle*{2}}
\put(56,20){\circle*{2}}
\put(11,5){\line(1,0){15}}
\put(41,5){\line(1,0){15}}
\put(56,5){\line(1,0){15}}
\put(56,5){\line(0,1){15}}
\put(10,0){$1$}
\put(25,0){$2$}
\put(37,0){$n-2$}
\put(52,0){$n-1$}
\put(55,22){$n$} 
\put(67,0){$n+1$}
\multiput(27,5)(2,0){7}{\line(1,0){1}}
\put(0,20){$D_{n+1}$:}
\end{picture}
\end{center}

\begin{center}
\setlength{\unitlength}{1mm}
\begin{picture}(86,24)
\put(11,12){\circle*{2}}
\put(26,12){\circle*{2}}
\put(41,12){\circle*{2}}
\put(56,12){\circle*{2}}
\put(71,12){\circle*{2}}
\put(11,12){\line(1,0){15}}
\put(56,12){\line(1,0){15}}
\put(10,7){$1$}
\put(25,7){$2$}
\put(40,7){$n$}
\put(51,7){$2n-2$}
\put(66,7){$2n-1$}
\multiput(27,12)(2,0){7}{\line(1,0){1}}
\multiput(42,12)(2,0){7}{\line(1,0){1}}
\put(0,20){$A_{2n-1}$:}
\end{picture}
\end{center}

\begin{center}
\setlength{\unitlength}{1mm}
\begin{picture}(86,24)
\put(11,5){\circle*{2}}
\put(26,5){\circle*{2}}
\put(41,5){\circle*{2}}
\put(56,5){\circle*{2}}
\put(71,5){\circle*{2}}
\put(41,20){\circle*{2}}
\put(11,5){\line(1,0){15}}
\put(26,5){\line(1,0){15}}
\put(41,5){\line(1,0){15}}
\put(56,5){\line(1,0){15}}
\put(41,5){\line(0,1){15}}
\put(10,0){$1$}
\put(40,22){$2$}
\put(25,0){$3$}
\put(40,0){$4$}
\put(55,0){$5$}
\put(70,0){$6$}
\put(0,20){ $E_6$:}
\end{picture}
\end{center}

\begin{center}
\setlength{\unitlength}{1mm}
\begin{picture}(86,24)
\put(26,5){\circle*{2}}
\put(41,5){\circle*{2}}
\put(56,5){\circle*{2}}
\put(41,20){\circle*{2}}
\put(26,5){\line(1,0){15}}
\put(41,5){\line(1,0){15}}
\put(41,5){\line(0,1){15}}
\put(25,0){$1$}
\put(40,22){$3$}
\put(40,0){$2$}
\put(55,0){$4$}
\put(0,20){$D_4$:}
\end{picture}
\end{center}
Their automorphism groups $\Gamma$ are given by the following generators $\bar{\sigma}$:
\begin{align}
\mbox{$D_{n+1}$: }\ &\bar{\sigma}(i)=i\, (1\le i \le n-1),\ \bar{\sigma}(n)={n+1},\ \bar{\sigma}({n+1})=n; \label{jw-12-21} \\ 
\mbox{$A_{2n-1}$: }\  &\bar{\sigma}(i)={2n-i}\, (1 \le i \le 2n-1);\label{jw-12-22}\\
\mbox{$E_6$: } \ &\bar{\sigma}(1)=6,\ \bar{\sigma}(3)=5,\ \bar{\sigma}(5)=3,
\ \bar{\sigma}(6)=1, \ \bar{\sigma}(j)=j\, (j=2, 4); \label{jw-12-23}\\
\mbox{$D_4$: }\  &\bar{\sigma}(1)=3,\ \bar{\sigma}(2)=2, \ \bar{\sigma}(3)=4, \ \bar{\sigma}(4)=1.\label{jw-12-24}
\end{align}
One can fold the Dynkin diagrams by using the above automorphisms to obtain the Dynkin diagrams of $B_n, C_n, F_4, G_2$ types.

Our main idea is that one should study integrable hierarchy mirror symmetry with a finite symmetry. First, we show that the above symmetry
$\Gamma$ can be endowed in both FJRW-theory and the Drinfeld-Sokolov hierarchies. The main result of the present paper is the following theorem.

\begin{thm}\label{liu-21}
The total descendant potential of  the $\Gamma$-invariant sector of  $D^T_{n+1}, A_{2n-1}, E_6$ FJRW-theory with the maximal diagonal symmetry group is a tau function of the corresponding $B_n, C_n, F_4$ Drinfeld-Sokolov hierarchy.
The total descendant potential of  the $\ZZ/3\ZZ$-invariant sector of $(D_4, \langle J\rangle)$ FJRW-theory  is a tau function of the  $G_2$ Drinfeld-Sokolov hierarchy.
\end{thm}

A key technical result is the following $\Gamma$-reduction theorem, which is of independent interest.
\begin{thm}\label{jw-1}
The $\Gamma$-invariant flows of an ADE Drinfeld-Sokolov hierarchy define the corresponding $B_n, C_n, F_4, G_2$ Drinfeld-Sokolov hierarchy. Furthermore, the
restriction of the ADE tau function to the $\Gamma$-invariant subspace of the big phase space provides a tau function of the corresponding $B_n, C_n, F_4, G_2$ Drinfeld-Sokolov hierarchy.
\end{thm}

We remark that the genus zero part of the $\Gamma$-invariant sector of FJRW-theory agrees with the $B_n, C_n, F_4, G_2$ Saito theory while the higher genus part is different from
the corresponding Givental-Dubrovin-Zhang theory. This explains the failure of $B_n, C_n, F_4, G_2$ Givental-Dubrovin-Zhang higher genus functions to yield tau functions of the corresponding Drinfeld-Sokolov integrable hierarchies.
Beyond the $B_n, C_n, F_4, G_2$ series, there are more Drinfeld-Sokolov hierarchies corresponding to twisted affine Kac-Moody algebras.
The integrable hierarchy mirror symmetry for the twisted series is still unresolved. We hope to come back to it on a different occasion.

The paper is organized as follows. In Section\,\ref{jw-12-3}, we will discuss the general set-up of cohomological field theories with a symmetry. In Section\,\ref{jw-12-4},
we will show how to endow the additional symmetry in FJRW-theory. In Section\,\ref{liu-4}, we present our proof of Theorem \ref{jw-3}, \ref{liu-21}, \ref{jw-1}. In Section\,\ref{jw-12-1}, we work out several explicit examples in which integrable
hierarchies are used to calculate FJRW-invariants. Our proof uses Fan-Jarvis-Ruan's theorem on the ADE integrable hierarchies conjecture. Fan-Jarvis-Ruan's theorem relies on early work of Frenkel-Givental-Milanov \cite{FGM, GM}.
We need to check that the chain of isomorphisms in the proof preserves $\Gamma$-action. This is true. However, it can be confusing for readers. In the Appendix, we present Dubrovin-Zhang's
alternative proof of the ADE integrable hierarchies conjecture, where the preservation of $\Gamma$-symmetry is transparent. 
In particular, Dubrovin-Zhang's proof bypasses Givental's higher genus theory
and the Kac-Wakimoto hierarchies. In a way, it is more direct.

\section{Cohomological field theory with a finite symmetry}\label{jw-12-3}

The first appearance of integrable hierarchies in Gromov-Witten
theory is the KdV-hierarchy. Its geometric counterpart is the intersection theory
on
the moduli space of stable Riemann surfaces $\overline{\Mcal}_{g,k}$.
The latter can be treated as the Gromov-Witten theory of the zero
dimensional manifold or FJRW-theory of $A_1=x^2$ singularity. Let's review it in more detail.

\subsection{Cohomological field theory}\label{sec-21}

Let $\overline{\Mcal}_{g,k}$ be  the moduli space of
isomorphism classes of genus $g$, stable, nodal Riemann
surfaces
with $k$ ordered markings. $\overline{\Mcal}_{g,k}$ is a central object in algebraic
geometry and has been studied intensively for decades. It is a smooth complex orbifold of
dimension $3g-3+k$. It is important to mention that $\overline{\Mcal}_{g,k}$ is only
well-defined in the so-called {\em stable range $2g+k\geq 3$}.
Over $\overline{\Mcal}_{g,k}$, each
marked point $x_i$ naturally defines an orbifold line bundle $L_i$
whose fiber at $C$ is $T^*_{x_i}C$. Let $\psi_i=c_1(L_i)$. One
can define the following intersection numbers for $\psi_i$ classes:
$$\langle \tau_{l_1}, \tau_{l_2}, \cdots,
\tau_{l_k}\rangle_g=\int_{\overline{\Mcal}_{g,k}}\prod_i \psi_i^{l_i} ,$$
defining them to be zero unless $\sum_i l_i=3g-3+k.$  They can be assembled into a generating function
$$\Fcal^g(t_0,t_1,\cdots)=\sum_{k\geq 0}\frac{t_{l_1}\cdots
t_{l_k}}{k!}\langle \tau_{l_1}, \tau_{l_2}, \cdots,
\tau_{l_k}\rangle_g,$$
which is a formal power series in infinitely many variables
$t_0, t_1, \cdots$. Then, we introduce {\em the total
descendant potential function}
$$\Dcal=\exp(\sum_{g\geq 0}\hbar^{g-1}\Fcal^g).$$
$\Dcal$ admits a geometric interpretation as the generating
function of intersection numbers for disconnected stable
Riemann surfaces.

We perform the dilaton shift
$$q_i=\left\{ \begin{array}{lll}
             t_i,&& i\neq 1;\\
             t_1-1,&& \mbox{otherwise.}
             \end{array}\right.$$
By the so-called dilaton equation, $\Fcal^g$ is a homogeneous power
series of degree $2-2g$ in the new variables $q_i$ for $g\ne 1$.
A central problem in mathematics and physics is to compute
$\Fcal^g$ or $\Dcal$. We can try to write them as combinations of known functions such as
exponential or trigonometric functions, or more general infinite products
such as modular forms or hypergeometric functions. If this
happens, we say that {\em $\Fcal^g$ or $\Dcal$ has a closed
formula}. Unfortunately, this almost never happens for
Gromov-Witten theory. The next-best situation is to find the
differential equations which it satisfies. We hope to find
enough equations to uniquely determine $\Fcal^g$ or $\Dcal$. Ideally, these equations are determined by the classical geometry
of the problem. It would be more striking if they came from entirely different sources.
The celebrated Witten-Kontsevich theorem is one such
example.

\begin{thm}
$\Dcal$ is a tau function of the KdV-hierarchy.
\end{thm}

\begin{rem}
$\Dcal$ is uniquely determined by the KdV-hierarchy
together with the so-called {\em string equation}.
\end{rem}

Since the KdV hierarchy is just the first example of a family of
hierarchies, it is natural to ask about the underlying geometry
of other integrable hierarchies. The basic idea is to consider more
general intersection numbers $\int_D\prod_i \psi^{l_1}_i$ for a
cycle $D\in H_*(\overline{\Mcal}_{g, k}, \QQ).$ This leads to
the interesting subject of special cycles of
$\overline{\Mcal}_{g,k}$. Furthermore, we require $D$ to
satisfy some general properties which are captured by the
notion of a cohomological field theory \cite{KM}.

Recall that there are several canonical morphisms between the
$\overline{\Mcal}_{g,k}$.
\begin{description}
\item[Forgetful Morphism]
$$\pi:
\overline{\Mcal}_{g,k+1}\rightarrow \overline{\Mcal}_{g,k}$$
forgets the last marked point $x_{k+1}$. Here, we
assume that $2g+k\geq 3$.
Furthermore, $\pi$ is the universal curve.
\item[Gluing the trees]
$$\rho_{tree}: \overline{\Mcal}_{g_1, k_1+1}\times
\overline{\Mcal}_{g_2, k_2+1}\rightarrow
\overline{\Mcal}_{g_1+g_2, k_1+k_2}.$$
\item[Gluing the loop]
$$\rho_{loop}: \overline{\Mcal}_{g, k+2}\rightarrow
\overline{\Mcal}_{g+1, k}.$$
\end{description}

Suppose that $H$ is a graded vector space with a nondegenerate
pairing $\langle\ , \ \rangle$ and a degree zero unit $1$. To simplify the
signs, we assume that $H$ has only even degree elements and the pairing
is symmetric. When $H$ has odd degree elements, everything
becomes ``super" and we leave it to readers to make the obvious
modifications.  Once and for all, we choose a homogeneous basis
$\phi_{\alpha}$ ($\alpha=1, \cdots, \dim H$) of $H$ with $\phi_1=1$. Let
$\eta_{\mu\nu}=\langle \phi_{\mu}, \phi_{\nu}\rangle$ and
$(\eta^{\mu\nu})=(\eta_{\mu\nu})^{-1}$.
\begin{defn}
A cohomological field theory is a collection of homomorphisms
$$\Lambda_{g,k}: H^{\otimes k}\rightarrow
H^*(\overline{\Mcal}_{g,k}, \QQ)$$
satisfying the following properties:

\begin{enumerate}
\item[\bf C0.] $\Lambda$ is homogeneous of degree
$$|\Lambda_{g,k}(a_1, \dots, a_k)|=N_{g,k}+\sum_i |a_i|.$$
for some rational numbers $N_{g,k}$.
\item[\bf C1.] The element $\Lambda_{g,k}$ is invariant under the
  action of the symmetric group $S_k$.

\item[\bf C2.] Let $g=g_1+g_2$ and $k=k_1+k_2$.  Then $\Lambda_{g,n}$ satisfies the composition property
\begin{multline} \label{eq:cfttree}
\rho_{\mathrm{tree}}^*
\Lambda_{g,k}(a_1, a_2,\ldots, a_k)=
\\
\Lambda_{g_1,k_1+1}(a_{i_1}, \ldots, a_{i_{k_1}},\phi_{\mu})\,
\eta^{\mu\nu} \otimes \Lambda_{g_2,k_2+1}(\phi_{\nu}
, a_{i_{{k_1}+1}},\ldots, a_{i_{k_1+k_2}})
\end{multline}
for all $a_i \in H$.
\item[\bf C3.] Let
\begin{equation}
\rho_{loop}: \overline{\Mcal}_{g-1, k+2}\rightarrow
\overline{\Mcal}_{g, k}
\end{equation}
be the gluing loop morphism. Then
\begin{equation}
\label{eq:cftloop}
\rho_{\mathrm{loop}}^*\,\Lambda_{g,k}(a_1, a_2,\ldots, a_k)\,=\,
\Lambda_{g-1,k+2}\,(a_1, a_2, \ldots, a_n, \phi_{\mu},
\phi_{\nu})\,\eta^{\mu\nu}
\end{equation}
for all $a_i \in H$.

\item[\bf C4a.] For all $a_i$ in $H$ we have
\begin{equation}
  \label{eq:identity}
\Lambda_{g,k+1}(a_1, \ldots, a_k, 1) =
\pi^*\Lambda_{g,k}(a_1,\ldots, a_k),
\end{equation}
where $\pi:\overline{\Mcal}_{g,n+1} \rightarrow
\overline{\Mcal}_{g,n}$ is the forgetful morphism. \item[\bf C4b.]
\begin{equation}
\label{eq:identity2}
\int_{\overline{\Mcal}_{0,3}}\,\Lambda_{0,3}(a_1, a_2,1)
= \langle a_1, a_2\rangle.
\end{equation}
\end{enumerate}
\end{defn}

\begin{rem}
The definition of cohomological field theory here is more restrictive than that of the original definition in the sense that  we require an additional  flat identity (C4a, C4b).
In the case of Gromov-Witten theory, we allow an additional
parameter $\beta$ to parametrize the curve class. Our main interest in this paper is FJRW-theory, where there is no curve class.
We leave the interested reader to modify the definition to suit the situation of Gromov-Witten theory.
\end{rem}

For each cohomological field theory, we can generalize the notion
of intersection numbers, the generating function and total
descendant potential function. Let
$$\langle \tau_{\alpha_1, l_1}, \cdots, \tau_{
\alpha_k, l_k}\rangle^{\Lambda}_g=\int_{\overline{\Mcal}_{g,k}}\prod_i \psi^{l_i}_i\Lambda_{g,k}(\phi_{\alpha_1},
\cdots, \phi_{\alpha_k}).$$
Associating a formal variable $t^{\alpha,p}$ to $\tau_{\al,p}$, we define the generating functions
$$\Fcal_{\Lambda}^g=\sum_{k\geq 0}\frac{t^{\alpha_1, l_1}\cdots
t^{\alpha_k, l_k}}{k!}\langle\tau_{\alpha_1, l_1}, \cdots,
\tau_{\alpha_k, l_k}\rangle_g^{\Lambda}$$
and its total descendant potential function
$$\Dcal_{\Lambda}=\exp(\sum_{g\geq
0}\hbar^{g-1}\Fcal^g_{\Lambda}).$$
Similarly, we can perform the dilaton shift
$$q^{\alpha, l}=\left\{\begin{array}{lll}
t^{\alpha, l},&& (\alpha, l)\neq (1,1);\\
t^{1, 1}-1,&& \mbox{otherwise}.
\end{array}\right.$$
$\Dcal_{\Lambda}, \Fcal_{\Lambda}$ can be viewed as formal functions on the {\em big phase space} $H^\infty=\mathrm{Spec}(\CC[[t^{\alpha,p}\mid \alpha=1, \dots, n;\ p=0, 1, 2, \dots]])$, which is the Cartesian
product of infinitely many copies of $H$:
\[H^\infty=H^{(0)}\times H^{(1)}\times H^{(2)}\times\cdots,\]
where $H^{(p)}=\mathrm{Spec}(\CC[[t^{\alpha,p}\mid \alpha=1, \dots, n]])$.

The goal is to find other cohomological field theories whose total descendant potential functions are
tau functions
of the given integrable hierarchies.

Consider the genus zero primary potential
$$F(v)=\Fcal^0_{\Lambda}|_{t^{\al,0}=v^\al,\ t^{\alpha, l}=0,\ l>0}$$
as a formal function in a neighborhood of $0\in H$. We trivialize the tangent bundle $TH$ near zero and view the pairing 
as a flat metric. 
One can  define a Frobenius algebra structure on $T_v H$ by the formula
$$\langle \frac{\p}{\p v^\al}\circ_v\frac{\p}{\p v^\beta}, \frac{\p}{\p v^\gamma}\rangle=\frac{\partial^3 F(v)}{\partial v^\al \partial v^\beta \partial v^\gamma}.$$
The associativity or WDVV equation follows from the gluing trees axiom.
The family of Frobenius algebra structures with a compatible flat metric yield a Frobenius manifold structure on $H$. The Frobenius manifold is generically semisimple if for a generic $v$ there is a basis $a_i$ such that 
$$a_i\circ_v a_j=\delta_{ij} a_i .$$
The examples considered in this article are generically semisimple.

Another important property is the Virasoro constraints (see the Appendix for the definition). The following theorem of Dubrovin-Zhang illustrates the importance of Virasoro constraints.
\begin{thm}
Suppose that $\Lambda=\{\Lambda_{g,k}\}$ is a generically semisimple cohomological field theory and satisfies the Virasoro constraints.  Then the higher genus 
potentials $\Fcal^g_\Lambda$ are uniquely determined by $F$ or the genus zero Frobenius manifold structure.
\end{thm}

\subsection{Cohomological field theory with a finite symmetry}

Suppose that a finite group $\Gamma$ acts on $H$ and preserves the grading, unit $1$ and pairing. A cohomological field theory with symmetry $\Gamma$
is a cohomological field theory 
$$\Lambda_{g,k}: H^{\otimes k}\rightarrow
H^*(\overline{\Mcal}_{g,k}, \QQ)$$
invariant under $\Gamma$  with trivial action on $H^*(\overline{\Mcal}_{g,k}, \QQ)$.
More generally, we can also allow a nontrivial action of $\Gamma$ on $H^*(\overline{\Mcal}_{g,k}, \QQ)$.

Let $H^{\Gamma}$ be the invariant subspace, and let
$$\Lambda_{g,k,\Gamma}: (H^{\Gamma})^{\otimes k}\rightarrow H^*(\overline{\Mcal}_{g,k}, \QQ)$$
be the restriction of $\Lambda_{g,k}$. Our key observation is the following:

\begin{pro}
$\Lambda_\Gamma=\{\Lambda_{g,k,\Gamma}\}$ satisfies all the axioms of cohomological field theory except the gluing loop axiom. Therefore, its genus zero theory
defines a Frobenius manifold structure on $H^{\Gamma}$ in a formal neighborhood of zero.
\end{pro}

\begin{prf} Only the gluing trees axiom needs a proof. We decompose $H=H^{\Gamma}\oplus H'$, where $H'$ is the direct sum of nontrivial irreducible factors.
The main trick is the following observation. Let $w\in H'$. Then, $ \sum_{g\in \Gamma} g(w)$ is $\Gamma$-invariant. Hence, it is in the intersection $H^{\Gamma}\cap H'=\{0\}$
and equals to zero. 
An easy consequence is that for $v\in H^{\Gamma}, w\in H'$, 
$$|\Gamma|\langle v, w\rangle=\sum_{\sigma} \langle \sigma (v), \sigma(w)\rangle=\sum_{\sigma} \langle v, \sigma(w)\rangle=\langle v, \sum_{\sigma} \sigma(w)\rangle=0.$$
Therefore, $H^{\Gamma}, H'$ are orthogonal with respect to the pairing. 

To prove the gluing trees axiom, suppose that $a_i \in H^{\Gamma}$.  Then we have
\begin{align*}
 &\rho_{\mathrm{tree}}^*\Lambda_{g_1+g_2,k}(a_1, a_2\,\ldots, a_k)\\
=&\Lambda_{g_1,k_1+1}(a_{i_1},\ldots, a_{i_{k_1}}, \phi_\mu)\eta^{\mu\nu} \otimes \Lambda_{g_2,k_2+1}(\phi_\nu
, a_{i_{{k_1}+1}},\ldots, a_{i_{k_1+k_2}})\\
&+
\Lambda_{g_1,k_1+1}(a_{i_1},\ldots, a_{i_{k_1}}, \phi_{\mu'})\eta^{\mu'\nu'} \otimes \Lambda_{g_2,k_2+1}(\phi_{\nu'}
, a_{i_{{k_1}+1}},\ldots, a_{i_{k_1+k_2}})
\end{align*}
for $\phi_\mu, \phi_\nu\in H^{\Gamma}, \phi_{\mu'}, \phi_{\nu'}\in H'$. Here we assume,
without lose of generality, that the basis $\phi_\alpha  (\alpha=1,\dots, \dim H)$ of $H$ is chosen w.r.t. the above decomposition of $H$. We claim that $\Lambda_{g,k}(a_1, \dots, a_i, w)=0$ for any $w\in H'$.  In fact, by using the same trick
as we used above, we get
$$|\Gamma|\Lambda_{g,k}(a_1, \dots, a_i, w)=\sum_{\sigma\in \Gamma} \Lambda_{g,k}(a_1, \dots, a_i, \sigma(w))$$
$$=\Lambda_{g,k}(a_1, \dots, a_i, \sum_{\sigma} \sigma(w))=0.$$
The proposition is proved.
\end{prf}

The above vanishing result shows that the existence of symmetry $\Gamma$ is a highly nontrivial property. 

\begin{defn}
If  a collection $\Lambda=\{\Lambda_{g, k}\}$ satisfies all the axioms except the gluing loop axiom, we call it {\em a partial} cohomological field theory.

\end{defn}
A partial cohomological field theory still defines a Frobenius manifold. From the above proposition, we know that $\Lambda_\Gamma=\{\Lambda_{g,k,\Gamma}\}$ is a partial cohomological field theory.
   
   The $\Gamma$-action on $H$ can be  lifted to the big phase space $H^\infty$.
Suppose the matrix of $\sigma\in \Gamma$ with respect to the basis $\{\phi_1, \dots, \phi_n\}$ is given by $A=(A^\beta_\alpha)$, i.e.
\[\sigma(\phi_\alpha)=A^\beta_\alpha \phi_\beta.\]
Then for any $f\in \Ocal_{H^\infty}$, we define the action of $\sigma$ on $f$ as follows:
\begin{equation}
\sigma^*(f)(t^{\alpha,p})=f((A^{-1})^\alpha_\beta t^{\beta,p}). \label{sigmastar}
\end{equation}
It is obvious that
\begin{lem}For a cohomological field theory with symmetry $\Gamma$, 
$\Dcal_{\Lambda}$ and $\Fcal_{\Lambda}$ are $\Gamma$-invariant.
\end{lem}

\section{FJRW-theory with a finite symmetry}\label{jw-12-4}

The main known examples of geometric cohomological field theories are the Gromov-Witten theory of a symplectic/K\"ahler orbifold $X$ 
and the FJRW-theory of an orbifolded singularity $(W,G)$. When $X$ admits a group action $\Gamma$, its Gromov-Witten cohomological field theory
has $\Gamma$ as its symmetry. 
However, the main examples connected to Drinfeld-Sokolov
hierarchies are FJRW-theory of ADE type.  The FJRW-theory of $(W,G)$ has a trivial $G$-action. It is not obvious how to endow a nontrivial
symmetry group $\Gamma$. In this section, we describe a method to enhance the symmetry of FJRW-theory by changing the
super potential.   A similar idea was also used in \cite{CIR}.   Let's first review FJRW-theory.

\subsection{Singularity $(W,G)$ and its FJRW state space}
\begin{defn}
$W$ is a \emph{quasi-homogeneous non-degenerate polynomial} if $W$ satisfies the following properties:
\begin{enumerate}
\item (Quasihomogeneity) There exist weights $q_i\in\mathbb{Q}$ such that for all $\lambda\in\mathbb{C}^*$,
\begin{equation*}
W(\lambda^{q_1}x_1,\cdots,\lambda^{q_N}x_N)=\lambda W(x_1,\cdots,x_N).
\end{equation*}
\item The choice of the weights $q_i$ is unique.
\item $W$ has an isolated singularity only at 0.
\end{enumerate}
\end{defn}
For each non-degenerate $W$, we have a symmetry group called \emph{maximal diagonal symmetry group} of $W$:
\begin{equation*}
G_{max,W}:=\Big\{\gamma=(\lambda_1,\cdots,\lambda_N)\in(\mathbb{C}^*)^N\Big|W(\lambda_1x_1,\cdots,\lambda_N x_N)=W(x_1,\cdots,x_N)\Big\}.
\end{equation*}
If  there is no confusion, we will often drop $W$ from the notation and denote it by $G_{max}$. $G_{max}$ always contains the \emph{exponential grading element} $J:=(e^{2\pi iq_1},\cdots,e^{2\pi iq_N})$.
In addition to $W$, FJRW theory depends on a choice of a suitable group $G$ with $\langle J\rangle \subset G \subset G_{max}$.

The \emph{central charge} of $W$ is defined to be
\begin{equation*}
\hat{c}_{W}:=\sum_{i=1}^{N}(1-2q_i).
\end{equation*}
For $\hat{c}_{W}<1$, $W$ is called a \emph{simple singularity}.  These have been completely classified into the famous ADE-singularities. 

For any $\gamma\in\,G$, we denote by $\mathbb{C}_{\gamma}^{N_{\gamma}}$ the fixed points of $\gamma$, where $N_{\gamma}$ is the complex dimension of the fixed locus. We denote by $W_{\gamma}$ the restriction of $W$ to the fixed locus. According to \cite{FJR1}, $W_{\gamma}$ is also non-degenerate.
\begin{defn}
We define the \emph{$\gamma$-twisted sectors} $H_{\gamma}$ to be the $G$-invariant part of the relative middle-dimensional cohomology for $W_{\gamma}$:
\begin{equation*}
H_{\gamma}:=H^*(\mathbb{C}^{N_{\gamma}}_{\gamma};W^{\infty}_{\gamma};\mathbb{C})^G.
\end{equation*}
Here $W^{\infty}_{\gamma}=(Re\,W_{\gamma})^{-1}(M,\infty)$ for $M\gg0.$

The \emph{FJRW state space} $H^*_{FJRW}(W,G)$ is defined to be the direct sum of all \emph{$\gamma$-twisted sectors $H_{\gamma}$} for the pair $(W,G)$:
\begin{equation}
H^*_{FJRW}(W,G):=\bigoplus_{\gamma\in G}H_{\gamma}=\bigoplus_{\gamma\in G} H^*(\mathbb{C}^{N_{\gamma}}_{\gamma};W^{\infty}_{\gamma};\mathbb{C})^G.
\end{equation}
\end{defn}
For any $\gamma\in G$, we have $\gamma=(\exp(2\pi i\Theta_1^{\gamma}),\cdots,\exp(2\pi i\Theta_N^{\gamma}))\in(\mathbb{C}^*)^N$ for some unique $\Theta_i^{\gamma}\in[0,1)\cap\mathbb{Q}$.
We define the \emph{degree shifting number} $\iota_{\gamma}$ by:
\begin{equation*}
\iota_{\gamma}:=\sum_{i=1}^N(\Theta_i^{\gamma}-q_i).
\end{equation*}
For each homogeneous element $\alpha\in H_{\gamma},$ the \emph{complex degree} $\deg_W(\alpha)$ is defined as
\begin{equation*}
\deg_W(\alpha):=\deg(\alpha)+\iota_{\gamma}.
\end{equation*}
Then $H^*_{FJRW}(W,G)$ is a graded vector space under this grading.
\begin{defn}
We say that the $\gamma$-twisted sector $H_\gamma$ is \emph{narrow} if $N_{\gamma}=0$. Otherwise, we say the $\gamma$-twisted sector is \emph{broad}.
\end{defn}
For each $\gamma\in\,G$, there is a natural intersection pairing $\langle\,,\rangle_{\gamma}$:
\begin{equation*}
\langle\,,\rangle_{\gamma}:H_{\gamma}\times\,H_{\gamma^{-1}}\longrightarrow\mathbb{C}.
\end{equation*}
Note that $H_{\gamma^{-1}}=H_{\gamma}$.
The pairing on narrow sectors is obvious since it is one-dimensional. 
A broad sector $H_\gamma$ with its pairing is isomorphic to ${\Omega}_{W_{\gamma}}$. 
Here for a quasi-homogeneous non-degenerate polynomial $W$ we define  
$\Omega_W:=\Omega^N/(d+d W)\Omega^{N-1}=\mathbb{C}[x_1,\cdots,x_N]/Jac(W)dx_1\dots dx_N$ with the residue pairing, where $Jac(W)$ is the Jacobi ideal of $W$. 
The residue pairing is given by
\begin{equation}\label{eq:Res}
\langle f,g\rangle:=\Res_{x=0}\frac{fg dx_1\cdots dx_N}{\frac{\partial W}{\partial x_{1}}\cdots
\frac{\partial W}{\partial {x_N}}}=C\mu,
\end{equation}
where $\mu$ is the dimension of $\Omega_W$  as a vector space and $C$ is the unique constant such that, modulo the Jacobi ideal, $fg=C\cdot Hessian(W)$.
The pairing on the FJRW state space is defined to be
\begin{equation*}
\langle\,,\rangle:=\sum_{\gamma\in\,G}\langle\,,\rangle_{\gamma}.
\end{equation*}

\subsection{W-Structure, virtual cycle and Cohomological Field Theory}
We start with a genus $g$ possibly nodal orbi-curve with $k$ marked points, denoted as  $\Sigma_{g,k}$, and its log canonical bundle
\begin{equation*}
\omega_{\log}:=\omega_{\Sigma_{g,k}}\bigotimes_{i=1}^{k}{\Ocal}(p_i)
\end{equation*}
for marked points $p_i$. We write the polynomial $W=\sum_{j=1}^{s}W_j$ as a sum of monomials $W_{j}=c_j\prod_{i=1}^{N}x_{i}^{b_{j,i}}$. A $W$-structure on the curve $\Sigma_{g,k}$ is a choice of $N$ orbifold line bundles ${\cal L}_1,\cdots,{\cal L}_N$ and $s$ isomorphisms of line bundles
\begin{equation*}
\varphi_j:\ W_j( {\cal L}_1,\cdots,{\cal L}_N)=\bigotimes_{i=1}^{N}{\cal L}_i^{\otimes\,b_{j,i}}\longrightarrow\omega_{\log}.
\end{equation*}
We denote the moduli space of $W$-structures as ${\Wcal}_{g,k}$. By \cite{FJR1}, the orbifold structure at a marked point
(or node) is specified by a group element $\gamma\in G_{max, W}$. By fixing the orbifold decorations at marked points, we have the decomposition:
\begin{equation*}
{\Wcal}_{g,k}=\sum_{\boldsymbol\gamma}{\Wcal}_{g,k}(\boldsymbol\gamma).
\end{equation*}
Here $\boldsymbol\gamma=(\gamma_1,\cdots,\gamma_k)$ is the orbifold decorations at marked points.
If we forget both the $W$-structure and the orbifold structure, then we have the forgetful morphism
\begin{equation*}
st:{\Wcal}_{g,k}\longrightarrow\overline{\Mcal}_{g,k}.
\end{equation*}
Here $\overline{\Mcal}_{g,k}$ is the Deligne-Mumford stack of stable curves.

The moduli of $W$-structures with group $G_{max}$  is a smooth Deligne-Mumford  stack denoted by ${\Wcal}_{g,k}$. In \cite{FJR1}, Fan-Jarvis-Ruan constructed a \emph{virtual fundamental cycle}:
\begin{equation*}
[{\Wcal}_{g,k}(\gamma)]^{vir}\in\,H_{*}({\Wcal}_{g,k},\mathbb{Q})\otimes\prod_{\tau\in\,T(\gamma)}H_{N_{\gamma_{\tau}}}(\mathbb{C}_{\gamma_{\tau}}^{N_{\gamma_{\tau}}},W_{\gamma_{\tau}}^{\infty},\mathbb{Q})^{G_{max}}
\end{equation*}
where $\gamma_{\tau}\in\,G_{max}$. Using the virtual cycle, they defined a cohomological field theory
$$\Lambda_{g,k}^{W,G_{max}}:(H^*_{FJRW})^{\otimes k}\longrightarrow H^*(\overline{\Mcal}_{g,k})$$
by
\begin{equation*}
\Lambda_{g,k}^{W,G_{max}}(\boldsymbol{a}):=\frac{|G_{max}|^g}{\deg(st)}PD\,st_*\Big([{\Wcal}_{g,k}(\gamma)]^{vir}\cap\prod_{i=1}^{k}a_i\Big)
\end{equation*}
for $a_i\in\,H_{\gamma_i}$, extended linearly to the whole space. Here, $PD$ is the Poincar\'e duality map.

For a so-called admissible subgroup $\langle J\rangle \subset G\subset G_{max}$, we can always choose a quasi-homogeneous Laurent polynomial  $Z$ of the same weights such that $G_{max, W+Z}=G$. Then, we apply the same construction for $W+Z$ to obtain $\Lambda_{g,k}^{W+Z,G}$. One can show that $\Lambda_{g,k}^{W+Z,G}$ is independent of $Z$ and hence can be denoted by $\Lambda_{g,k}^{W,G}$.

\begin{thm}[Fan--Jarvis--Ruan \cite{FJR2}]\label{thm:CohFT}
We denote the FJRW state space by $\Hcal_{W,G}$.
Let $1$ be the distinguished generator $1_{J}$ attached to the exponential grading element $J=exp(2\pi i q_1, \cdots, 2\pi i q_N)$  lying in the $A$-admissible group $G$ (i.e. $J\in G$). Let
$\langle\cdot ,\cdot\rangle^{W,G}$ denote the pairing on $\Hcal_{W,G}$. Then, the collection 
$$(\Hcal_{W, G}, \langle\cdot ,\cdot \rangle^{W,G},
\{\Lambda^{W,G}_{g,k}\}, 1)$$ is a cohomological field theory.
\end{thm}

The following properties hold true:
\begin{enumerate}
                                \item {Decomposition.}
If $W_1$ and $W_2$ are two singularities in distinct
variables, then the cohomological field theory arising from $(W_1+W_2,G_1\times G_2)$ is
the tensor product of the cohomological field theories arising
from $(W_1,G_1)$ and $(W_2,G_2)$.

\item {Deformation invariance.}
Suppose that $W_t, t\in [0,1]$ is a one-parameter family of nondegenerate polynomials such that $W_t$ is
$G$-invariant. Then we have a
canonical isomorphism  $\Hcal_{W_0, G}\cong \Hcal_{W_1, G}$.
Under the above isomorphism,
$$\Lambda^{W_0, G}_{g, k}=\Lambda^{W_1, G}_{g, k}.$$
Namely, $\Lambda^{W,G}_{g,k}$ depends only on $(q_1, \dots, q_N)$ and on $G$.
Note also that, when applied to a deformation of a polynomial $W$ along a loop,
this property implies monodromy invariance for  $\Lambda^{W,G}_{g,k}$.
\item{$G_{max}$-invariance.}
The group $G_{max}$ acts on $\Hcal_{W,G}$
in an obvious way.
Then, $\Lambda^{W,G}_{g,k}$ is invariant with respect to the action of
$\Aut(W)$ of diagonal symmetries on each state space entry $a_1,\dots, a_n\in \Hcal_{W,G}$.
\end{enumerate}

The present paper needs to use the following theorems. 

\begin{thm}
\begin{description}
\item[(1)] The genus zero FJRW-theories of $D^T_n, A_{2n-1},$ and $E_6$ potentials with $G_{max}$ are isomorphic to the Saito Frobenius manifold of the transpose singularities $W^T$.
The genus zero FJRW-theory of $(D_4, \langle J \rangle)$ is isomorphic to the Saito Frobenius manifold of the transpose singularities of $D_4$.
\item[(2)] The all genus  FJRW-theory of $ADE, D^T$-potential with $G_{max}$ and  ($D_{2n}, \langle J\rangle)$ satisfy the Virasoro constraints. 
\end{description}

\end{thm}

Property (2) of the above theorem is a consequence of Theorem 6.1.3 of \cite{FJR2} and the proof of Theorem 1.1 of \cite{FJMR}. Alternatively, one can use Teleman's solution of Givental conjecture \cite{Te}. Property (1) is their genus zero part.

\begin{rem}
The all genus FJRW-theory of the above examples are isomorphic to Givental's higher genus theory of the mirror singularities, but we do not need this fact in the current paper. Instead, we use
Dubrovin-Zhang's axiomatized integrable system theory to bypass the work of Frenkel-Givental-Milanov \cite{FGM, GM}.
\end{rem}

\subsection{Enhance the symmetry of FJRW-theory}\label{liu-3}

Let's set up the notation.   If $\gamma$ is narrow, $\Hcal_\gamma$ is one-dimensional and we label the canonical
generator by $e_{\gamma}$. For broad sectors, we use the $G$-equivariant isomorphism
$$H^{mid}(\CC_\gamma^{N_\gamma}, W_\gamma^{\infty}, \CC)^G\cong \Omega^G_{W_\gamma}$$
to label the generator by $e_{\gamma} \phi $ for $\phi\in \Omega_{W_\gamma}^G$.

\begin{exa}
Let's start from $A_{2n-1}=x^{2n}$. Then $G_{max}=\ZZ/2n\ZZ$, generated by $J=e^{\frac{\pi i}{n}}$. The state space $\Hcal_{A_{2n-1}}$ is generated by narrow
sectors $e_{J^i}$ for $0<i<2n$ with degrees given by  $\deg(e_{J^i})=\frac{i-1}{2n}$. We claim that $\sigma (e_{J^{i+1}})=(-1)^{i}e_{J^{i+1}}$ is a symmetry of the cohomological field theory.
Note that 
$$\Lambda_{g,k}^{A_{2n-1}, G_{max}}(\sigma( e_{J^{i_1+1}}), \cdots, \sigma(e_{J^{i_k+1}}))=(-1)^{\sum_j i_j}\Lambda_{g,k}^{A_{2n-1}, G_{max}}( e_{J^{i_1+1}}, \cdots, e_{J^{i_k+1}}).$$
We claim that 
$\Lambda_{g,k}^{A_{2n-1}, G_{max}}( e_{J^{i_1+1}}, \cdots, e_{J^{i_k+1}})=0$ if $\sum_j i_j$ is odd and hence $\sigma$ preserves $\Lambda_{g,k}^{A_{2n-1}, G_{max}}$.
To prove the claim, we use the following geometric  property of $\Wcal_{g,k}(J_{i_1+1}, \cdots, J_{i_k+1})$. Suppose that $L^{2n}\cong \omega_{log}$ is a $2n$-spin structure. Then, 
$$\deg (|L|)=\frac{2g-2-\sum_j i_j}{2n}\in \ZZ.$$
In particular, $\Wcal_{g,k}(J_{i_1+1}, \cdots, J_{i_k+1})$ is empty if  $\sum_j i_j$ is odd.
\end{exa}

\begin{exa}
Another example of similar nature is $E_6=x^3+y^4$. $G_{max}$ is generated by $J=(\eta^4, \eta^3)$ for $\eta^{12}=1$.
The state space is generated by the narrow sectors $e_{J}, e_{J^2}, e_{J^5}, e_{J^7}, e_{J^{10}}, e_{J^{11}}$. Consider the action
$$\sigma: e_{\gamma}\rightarrow (-1)^{4\Theta_y^{\gamma}-1} e_{\gamma}.$$
More precisely, we multiply by $-1$ on $e_{J^2},e_{J^{10}}$ and act trivially otherwise.
Then we have
$$\Lambda_{g,k}^{E_6, G_{max}}(\sigma( e_{\gamma_1}), \cdots, \sigma(e_{\gamma_k}))=(-1)^{\sum_i (4\Theta^{\gamma_i}_y-1)}\Lambda_{g,k}^{E_6, G_{max}}( e_{\gamma_1}, \cdots, e_{\gamma_k}).$$
In this case, an element of $\Wcal_{g,k}(\gamma_1, \cdots, \gamma_k)$ is a pair of orbifold line bundles $L_x, L_y$ such that $L^3_x\cong L^4_y\cong \omega_{log}$, satisfying some
additional constraints.  The fact that $\deg( |L_y|)\in \ZZ$ implies that $\sum_i (4\Theta^{\gamma_i}_y-1)$ has to be even. 
\end{exa}

The previous cases follow from degree considerations. For other cases, such a simple argument fails. The main new idea is as follows.
By the deformation invariance axiom, the FJRW-theory of $(W,G)$ depends only on the weights $q_i$ and $G$. In particular, it is independent of the choice of $W$.
On the other hand, $G_{max}$ does depend on $W$. This provides an interesting possibility that we can choose a different $\tilde{W}$ with the same weights such that 
its $G_{max, \tilde{W}}$ is bigger than $G_{max, W}$.  Then the $G_{max, \tilde{W}}$-invariance axiom 
tells us that our theory is in fact $G_{max, \tilde{W}}$-invariant.  It follows that $\Gamma=G_{max, \tilde{W}}/G_{max, W}$ is a symmetry of the FJRW-theory of $(W,G)$! Let's put this idea into practice.

\begin{exa}
Consider $D^T_{n+1}=x^ny+y^2$. Its weights are $q_x=\frac{1}{2n}, q_y=\frac{1}{2}$. $G_{max}=\ZZ/2n\ZZ$ is generated by $J=(\xi, -1)$ for $\xi^{2n}=1$.
On the other hand, if we choose $\tilde{D}^T_{n+1}=x^{2n}+y^2$, it has the same weights but $\tilde{G}_{max}=\ZZ/2n \ZZ\times \ZZ/2\ZZ$. Therefore, 
$(D^T_{n+1}, G_{max})$ has an additional $\ZZ/2\ZZ$ symmetry.

The FJRW state space $\Hcal_{D^T_{n+1}, G_{max}}$ consists of 
narrow sectors $\Hcal_{J^i}$ for each odd integer $i\leq 2n-1$ and a single 
broad sector $ \Hcal_{J^0}$.
Again, the extra $\ZZ/2\ZZ$ symmetry acts trivially on narrow sectors and as multiplication by $\pm 1$ on the broad sector.

\end{exa}

\begin{exa}
Consider $D_4=x^3+xy^2$. Its weights are $q_x=\frac{1}{3}, q_y=\frac{1}{3}$ and $\langle J \rangle=\ZZ/3\ZZ$. On the other hand, we can
choose $\tilde{D}_4=x^3+y^3$. The new $\tilde{G}_{max}=\ZZ/3\ZZ\times \ZZ/3\ZZ$. Therefore, the FJRW-theory of $(D_4, \langle J \rangle)$ has an additional symmetry $\Gamma=\tilde{G}_{max}/\langle J \rangle=\ZZ/3\ZZ$.

The FJRW state space $\Hcal_{D_4, \langle J \rangle}$ consists of the narrow sectors $e_{J}, e_{J^2}$ and a broad sector
$\Hcal_{J^0}=\Omega_{D_4}^{\langle J \rangle}$ of dimension two generated by $ e_{J^0}x dxdy$ and $e_{J^0}ydxdy$. The extra $\ZZ/3\ZZ$ symmetry acts trivially on the narrow sectors and it acts on the broad sector by 
$$(e_{J^0}xdxdy, e_{J^0}ydxdy)\rightarrow (\xi e_{J^0}xdxdy, \xi^{-1}e_{J^0}ydxdy ).$$
Here $\xi^3=1$. In \cite{FJMR}, it is proved that the genus zero 
FJRW-theory of $(D_4, \langle J \rangle)$ is isomorphic to the Saito-theory of $D_4$. 

\end{exa}

We summarize the results of this section in the following theorem.

\begin{thm}
The FJRW-theory of $D^T_{n+1}, A_{2n-1}, E_6$  with $G_{max}$ and $(D_4, \langle J\rangle)$ have the following nontrivial $\Gamma$ symmetries: (1) $\Gamma=\ZZ/2\ZZ$ for $A_{2n-1}, E_6, D^T_{n+1}$ with $G_{max}$; (2) $\Gamma=\ZZ/3\ZZ$ for 
$(D_4, \langle J \rangle)$.
\end{thm}

\section{Proof of the Main Theorem}\label{liu-4}

The idea of the proof is as follows. From the last section, $\Dcal_{\Lambda}$ and $\Fcal_{\Lambda}$ (as functions on the big phase spaces) and the associated Frobenius manifolds for the FJRW-theory of $D^T_n, A_{2n-1}, E_6$ 
 with $G_{max}$ and $(D_4, \langle J \rangle)$ are $\Gamma$-invariant. 
 We view FJRW-theory as the $A$-model. Next, we work on the Drinfeld-Sokolov hierarchies, which we consider as the $B$-model.  In the B-model setting,
$\Gamma$ originates from symmetries of the Dynkin diagram. 
In \cite{DS}, Drinfeld and Sokolov constructed an integrable hierarchy from any affine Lie algebra $\hg$ and a chosen vertex of
its Dynkin diagram. When $\hg$ is untwisted and the chosen vertex is the zeroth one, the Drinfeld-Sokolov hierarchy has
a semisimple bihamiltonian structure of hydrodynamic type whose leading term is given by the Frobenius manifold structure defined
on the orbit space of the corresponding Coxeter group \cite{DLZ2}. When the affine Lie algebra $\hat{\g}$ is of ADE type, the Frobenius manifold structure coincides with the one that is defined on the space of miniversal deformations of the corresponding simple singularity of ADE type \cite{S1}. 

The symmetry $\Gamma$ of the Dynkin diagram induces an action on the associated Drinfeld-Sokolov hierarchy. The action on the Drinfeld-Sokolov hierarchy then induces an action on its phase space and on the Frobenius manifold 
corresponding to the Coxeter group/miniversal deformation of ADE singularities.
We can calculate explicitly the induced $\Gamma$-action on the Frobenius manifold associated to the Coxeter group/miniversal deformation of ADE singularties.
Then, it is easy to use the mirror theorem of Fan-Jarvis-Ruan to match the $\Gamma$-action on FJRW-theory with the action on the mirror Saito Frobenius manifold induced from Drinfeld-Sokolov hierarchy. 
Hence, the total descendant potential function $\Dcal_{\Lambda}$ of FJRW-theory
can be viewed as a $\Gamma$-invariant function on the mirror B-model big phase space. By the theorem of Fan-Jarvis-Ruan, it is a ($\Gamma$-invariant) tau function of the mirror Drinfeld-Sokolov hierarchy.
The key technical work of our paper is the $\Gamma$-reduction theorem where we show that the invariant flows of ADE Drinfeld-Sokolov hierarchies define the BCFG Drinfeld-Sokolov hierarchies. 
Moreover, the restriction of  the tau functions to the  $\Gamma$-invariant subspaces provide the tau functions of BCFG hierarchies.

In this section, we first review the Drinfeld-Sokolov hierarchies with an emphasis on their bihamiltonian structures. In section 4.1.2, we review Dubrovin-Liu-Zhang's proof that Givental's total descendant potential functions of the BCFG Frobenius manifolds are not tau functions of the corresponding Drinfeld-Sokolov hierarchies, which directly motivated our current work. In Section 4.2,
we prove the $\Gamma$-reduction theorem. In Section 4.3, we calculate the $\Gamma$-action on the Frobenius manifold, corresponding to the Coxeter group/miniversal deformation of ADE singularities, that is induced from the $\Gamma$-action on the associated Drinfeld-Sokolov hierarchy and tie all the arguments together.

\subsection{Review of Drinfeld-Sokolov hierarchies}\label{sec-4-1}

In this paper, we only consider Drinfeld-Sokolov hierarchies that are associated to untwisted affine Lie algebras with the zeroth
vertex of their Dynkin diagrams being chosen. Our construction is mainly based on that of \cite{DS, Wu2}.

Let $\g$ be a simple Lie algebra of rank $n$, and let $(\ \mid\ )$ be the normalized Killing form of $\g$ such that the longest roots of
$\g$ have square length $2$. Let the loop algebra $L(\g)$ be given by $L(\g)=\CC[\lambda, \lambda^{-1}]\otimes \g$.  Then the algebra
$\hg$ is defined as the universal central extension 
\[0\to \CC K \to \hg \to L(\g)\to 0\]
of $L(\g)$ with the Lie bracket
\[[X(\lambda)+a\,K, Y(\lambda)+b\,K]_{\hg}=[X(\lambda),Y(\lambda)]_{L(g)}+\Res_{\lambda=0}(X'(\lambda)\mid Y(\lambda))K.\]
Note that $\hg$ does not contain the derivation $d=\lambda \frac{d}{d\lambda}$ here (c.f. \S 7.2 of \cite{Kac}).

Let $\h$ be a Cartan subalgebra of $\g$.  $\g$ has the root space decomposition
\[\g=\h\oplus \bigoplus_{\alpha\in\Delta} \g_{\alpha},\]
where $\Delta\subset \h^*$ is the root system of $\g$. Fix a basis $\Pi=\{\alpha_1, \dots, \alpha_n\}$ of $\Delta$.
One can choose a set of Weyl generators $E_i\in \g_{\alpha_i}$, $F_i\in \g_{-\alpha_i}$, and
$H_i=[E_i, F_i]$ such that
$(E_i | F_i)=2/(\alpha_i | \alpha_i)$. Let $\theta\in\Delta$ be the highest root with respect to $\Pi$.
We choose $F_0\in\g_{\theta}$, $E_0\in\g_{-\theta}$ such that
\[(E_0 | F_0)=2/(\theta | \theta),\quad F_0=-\omega(E_0),\]
where $\omega$ is the Chevalley involution defined by the following relations:
\[\omega(E_i)=-F_i,\quad \omega(H_i)=-H_i, \quad \omega(F_i)=-E_i.\]
The Lie algebra $\hg$ is generated by 
\begin{align*}
& e_0=\lambda\otimes E_0, \quad f_0=\lambda^{-1}\otimes F_0, \\
& e_i=1\otimes E_i, \quad f_i=1\otimes F_i,\quad i=1, \dots, n,
\end{align*}
and the central element $K$.

For an arbitrary $s=(s_0, s_1, \dots, s_n) \in \ZZ^{n+1}_{\ge 0}$, one can define a gradation on $\hg$ by
\[\deg_s e_i=-\deg_s f_i=s_i\ (i=0, \dots, n),\quad \deg_s K=0.\]
In the present paper, the following two gradations are frequently used:
\begin{itemize}
\item[i)] the homogeneous gradation: $s'=(1, 0, \dots, 0)$. We denote
\[\hg_j=\{X\in \hg\mid \deg_{s'}X=j\},\ j\in \ZZ.\]
\item[ii)] the principal gradation: $s''=(1, 1, \dots, 1)$.  We denote
\[\hg^j=\{X\in \hg\mid \deg_{s''}X=j\},\ j\in \ZZ.\]
\end{itemize}
We will also use notations such as $\hg_{\ge 0}=\bigoplus\limits_{j\ge 0}\hg_j$, $\hg^{< 0}=\bigoplus\limits_{j< 0}\hg^j$.
In particular, $\hg_0=\g\oplus \CC K$ and $\hg^0=\h\oplus\CC K$. Here, we regard $\g$ and $\h$ as the subalgebras
$1\otimes \g$ and $1\otimes \h$ of $\hg$.
Any $X\in \hg$ can be uniquely decomposed as $X=X_{+}+X_{-}$,
where $X_{+}\in \hg_{\ge 0}$ and $X_{-}\in \hg_{<0}$. We denote the projections $X\mapsto X_{+}$ and $X\mapsto X_{-}$
by $(\ )_+$ and $(\ )_-$ respectively.

Let
\begin{equation}\label{zh-Lambda}
\Lambda=\sum\limits_{i=0}^n e_i
\end{equation}
and consider the 
subalgebra $\s=\Ker \ad_\Lambda$ of $\hg$.
Let $E$ be the set of exponents of $\hg$ (see \S 14.3 of \cite{Kac} for their values). 
We can choose a basis $\{K\}\cup\{\Lambda_j\in \hg^j\}_{j\in E}$ of $\s$
such that 
\[[\Lambda_i, \Lambda_j]=i \delta_{i+j,0}K.\]
The subalgebra $\s$ is called the principal Heisenberg algebra associated to $\Lambda$.

\subsubsection{The Drinfeld-Sokolov hierarchies and their tau functions}
\label{sec-411}

Now let us take a matrix realization of $\g$, so that elements of $\g$ become complex matrices. Since
$\hg=L(\g)\oplus \CC K$, every element $X$ of $\hg$ can be uniquely represented as $X=X_1+ X_2$, where $X_1$ is
a matrix whose entries are Laurent polynomials in $\lambda$, and $X_2$ is a constant multiple of $K$. From now on,
the projection $X\mapsto X_2$ will be denoted by $(\ )_K$, and we will denote $X\in L(\g)$ if $(X)_K=0$.

Let us introduce a matrix differential operator
\begin{equation}
\Lcal=\frac{d}{d x}+\Lambda+q,\quad q\in C^\infty(\RR, \g\cap \hg^{\le 0}), \label{lax}
\end{equation}
where $x$ is the coordinate of $\RR$.

\begin{pro}[\cite{DS, Wu2}]\label{pro-51}
For the matrix differential operator $\Lcal$ given in \eqref{lax}, there exists a unique pair of matrices
\[U\in C^\infty(\RR, \hg^{<0}) \quad \mbox{and} \quad H\in C^\infty(\RR, \s\cap \hg^{<0})\]
such that
\begin{itemize}
\item[i)] $e^{-\ad_U}(\Lcal)=\frac{d}{dx}+\Lambda+H$.
\item[ii)] For any $j\in E_+$, $e^{\ad_U}(\Lambda_j)\in L(\g)$,
\end{itemize}
where $E_+$ is the set of positive exponents (see \S 14.3 of \cite{Kac}).
Moreover, the entries of the matrices $U$ and $H$ are differential polynomials of the entries of the matrix $q$.
\end{pro}

For each $j\in E_+$, one can show that the equation
\[\frac{\p \Lcal}{\p t_j}=-[\left(e^{\ad_U}(\Lambda_j)\right)_+, \Lcal]\]
defines a system of evolutionary PDEs for the entries of $q$, and moreover
\[[\frac{\p}{\p t_j}, \frac{\p}{\p t_k}]=0, \mbox{ for any } j, k \in E_+,\]
so the flows $\{\p/\p t_j\}_{j\in E_+}$ form an integrable hierarchy of evolutionary PDEs. We call this integrable hierarchy the {\em pre-Drinfeld-Sokolov hierarchy}.

We define the gauge transformations on the space of matrix differential operators of the form \eqref{lax} by $\Lcal\mapsto \tilde{\Lcal}=e^{\ad_W}(\Lcal)$, where 
$W\in C^\infty(\RR, \g\cap \hg^{<0})$. 
It is shown in \cite{DS} that the flows of the pre-Drinfeld-Sokolov hierarchy preserve gauge equivalence,  so they can be reduced to the space of gauge equivalence classes.
The reduced integrable hierarchy is called the {\em Drinfeld-Sokolov hierarchy}.

Let $\fb=\g\cap \hg^{\le 0}$ and $\fn=\g\cap \hg^{<0}$ be the Borel subalgebra and nilpotent subalgebra of $\g$
associated to $\h$ and $\Pi$,
and denote $I_+=\sum\limits_{i=1}^n E_i$. In \cite{DS}, Drinfeld and Sokolov showed that $\ad_{I_+}:\fn\to\fb$ is
injective. A subspace $V$ of $\fb$ such that $\fb=[I_+, \fn]\oplus V$ is called a gauge. It is easy to see that
the space of $q$ is just $C^\infty(\RR,\fb)$, and the space of gauge transformations is $C^\infty(\RR, \fn)$, so
the space of gauge equivalence classes can be identified with $C^\infty(\RR, V)$.

In this paper, we will use the lowest weight gauge, which is introduced in \cite{BF}. Let $a_{ij}=\alpha_j(H_i)$
be the $(i,j)$-th entry of the Cartan matrix of $\g$, and let $x_i$ be the solution to the equations
\[\sum_{i=1}^n x_i a_{ij}=1,\quad j=1, \dots, n.\]
Define
\[\rho=\sum_{i=1}^n x_i H_i,\quad I_-=\sum_{i=1}^n 2\,x_i\,F_i.\]
Then one can check that
\[[\rho, I_\pm]=\pm I_\pm,\quad [I_+, I_-]=2\rho,\]
so $\{I_+, 2\rho, I_-\}$ forms an $sl_2$ Lie algebra, and $\g$ becomes a module of this Lie algebra. Denote
\begin{equation}\label{fix-gauge}
V=\Ker \ad_{I_-},
\end{equation}
then it is easy to see that $V$ is a gauge. Let $1=m_1\le m_2 \le\dots\le m_n$ be the exponents of $\g$
(see e.g. \S 14.2 of \cite{Kac}).
One can show that the eigenvalues of $\rho$ on $V$ are exactly given by $\{-m_1, \dots, -m_n\}$. Let $\gamma_1, \dots, \gamma_n$
be a set of eigenvectors of $\rho$ corresponding to the above eigenvalues. Then the gauge fixed $\Lcal$ can be written as
\begin{equation}
\Lcal^{\can}=\frac{d}{dx}+\Lambda+q^{\can}=\frac{d}{dx}+\Lambda+\sum_{i=1}^n u^i(x)\gamma_i. \label{lax-can}
\end{equation}
Here the superscript ``$\can$'' stands for ``canonical'', and $u^i$ are gauge invariant differential polynomials of the entries of the matrix $q$ that appear in the 
operator $\Lcal$ given in  \eqref{lax}.
Note that each $\gamma_i$ generates an irreducible lowest weight module
\begin{equation}
L_i=\mathrm{Span}\{\gamma_i, \ad_{I_+}(\gamma_i), \dots,  \ad_{I_+}^{2m_i}(\gamma_i)\} \label{eq-Li}
\end{equation}
of the $sl_2$ Lie algebra $\{I_+, 2\rho, I_-\}$,
and $\g$ can be written as the direct sum of these irreducible submodules. This gauge is called the {\em lowest weight gauge}.
By reducing the pre-Drinfeld-Sokolov hierarchy to the space of $\Lcal^{\can}$, one obtains the Drinfeld-Sokolov hierarchy
with the lowest weight gauge, which is a hierarchy of  evolutionary PDEs for the unknown functions $u^1, \dots, u^n$.

Now let us consider the tau functions of the Drinfeld-Sokolov hierarchy.
\begin{pro}[\cite{Wu2}]\label{zh-12-2}
Let the matrix $H$  be given in Proposition \ref{pro-51} for $\Lcal$. Define
\[h_j=-j\frac{(\Lambda_j | H)}{(\Lambda_j | \Lambda_{-j})}, \quad j\in E_+.\]
Then $h_j$ can be represented as differential polynomials of $u^1, \dots, u^n$. They are Hamiltonian densities of the Drinfeld-Sokolov hierarchy.
Moreover, they satisfy the following tau symmetry properties:
\begin{equation}\label{jw-2}
\frac{\p h_j}{\p t_k}=\frac{\p h_k}{\p t_j}=\p_x\Omega_{jk},\quad j,k\in E_+,
\end{equation}
where $\Omega_{jk}$ are some differential polynomials of $u^1, \dots, u^n$.
\end{pro}

Define
\[\deg u^i=m_i+1,\quad \deg \p_x=1,\]
then we can uniquely determine the differential polynomials $\Omega_{jk}$ by the homogeneity condition $\deg \Omega_{jk}=m_j+m_k$ and the relations given in \eqref{jw-2}.
These differential polynomials are called the two-point functions of the Drinfeld-Sokolov hierarchy. They have
the following properties:
\[\Omega_{jk}=\Omega_{kj}, \quad \frac{\p \Omega_{jk}}{\p t_l}=\frac{\p \Omega_{jl}}{\p t_k}.\]
Suppose $u(x,t)=\{u^1(x,t), \dots, u^n(x,t)\}$ is a solution of the Drinfeld-Sokolov hierarchy, then the above conditions
imply the existence of a function
$\Fcal(x,t)$ such that
\begin{equation}\label{zh-12-1}
\frac{\p^2 \Fcal(x,t)}{\p t_j\p t_k}=\Omega_{jk}(u(x,t)).
\end{equation}
In \cite{DZ3}, the set of differential polynomials $\Omega_{i, j}, i, j\in E_+$ is called a \emph{tau-structure of the Drinfeld-Sokolov hierarchy}.
If $u(x,t)$ is the topological solution of the Drinfeld-Sokolov hierarchy that is specified by the string equation (see \cite{DZ3} and the Appendix for a detailed explanation), then the function $\Fcal$ is uniquely determined by the relation \eqref{zh-12-1} and the string equation up to the addition of a constant.
\begin{defn}
The function $\Fcal$ is called the free energy of $u(x,t)$,
and $\tau=e^\Fcal$ is called the tau function of $u(x,t)$.
\end{defn}

\subsubsection{The bihamiltonian structure}\label{jw-12-2}

The Drinfeld-Sokolov hierarchy constructed in the last subsection has the following bihamiltonian formalism:
\begin{equation}\label{jw-12-30}
\frac{\p u^i}{\p t_k}=\{u^i(x), c_{k,1} H_k\}_1=\{u^i(x), c_{k,2} H_{k+h}\}_2,
\quad i=1,\dots, n,\, k\in E_+.
\end{equation}
Here the densities $h_k$ of the Hamiltonians $H_k=\int h_k dx$ are defined in Proposition \ref{zh-12-2}, $h$ is the Coxeter number of the simple Lie algebra $\g$,
and $c_{k,a}\ (a=1, 2)$ are certain constants 
depending on the choice of the basis $\{\Lambda_j\}_{j\in E_+}$ of the principal Heisenberg subalgebra $\s$.
The brackets $\{\, \,,\, \}_a\ (a=1, 2)$ are two compatible Poisson brackets of partial differential equations.
They form the so-called Drinfeld-Sokolov bihamiltonian structure of the integrable hierarchy associated to $\g$.
For more details on Hamiltonian structures of PDEs and their compatibility, we refer readers to \cite{DZ3, LZ2}.

The bihamiltonian structure $\{\, \,,\, \}_{a}\ (a=1,2)$ has the following form:
\begin{align}
&\{u^i(x), u^j(y)\}_a=g^{ij}_a(u(x))\delta'(x-y)+\Gamma^{ij}_{k,a}(u(x))u^{k}_x(x)\delta(x-y)\nonumber\\
&\quad+\sum_{m\ge1}\sum_{k=0}^{m} P^{ij}_{m,k,a}(u(x),u_x(x),u_{xx}(x),\dots, \p_x^k u(x))\delta^{(m-k)}(x-y), \label{biham}
\end{align}
where $(g^{ij}_1), (g^{ij}_2)$ form a pair of compatible contravariant metrics (see \cite{Du2} for details), and $\Gamma^{ij}_{k,a}\ (a=1,2)$ are the contravariant Christoffel coefficients associated to the metrics.
The coefficients $P^{ij}_{m,k,a}\ (a=1,2)$ are certain homogeneous polynomials of $u^i_x(x), u^i_{xx}(x), \dots, \p_x^k u^i(x)$ satisfying the property
\[\deg P^{ij}_{m,k,a}=k, \mbox{ if we define } \deg \p_x^k u^i(x)=k.\]
The first two terms of the right-hand side of \eqref{biham} are called the leading terms of the bihamiltonian structure, while the second line of \eqref{biham} is
called the deformation part.

Denote the right-hand side of \eqref{jw-12-30} by $A_k^i(u, u_x, u_{xx}, \dots)$.
We have the dispersionless Drinfeld-Sokolov hierarchy consisting of the flows
\begin{equation}
\frac{\p v^i}{\p t_k}=\lim_{\hbar\to0}\frac{1}{\sqrt{\hbar}}A_k^i(v, \sqrt{\hbar} v_x, \hbar v_{xx}, \dots)=B^i_{kl}(v)v^l_x. \label{semi-classical}
\end{equation}
We can view the Drinfeld-Sokolov hierarchy as a certain deformation of the above integrable hierarchy of hydrodynamic type which has a bihamiltonian structure of the form \eqref{liu-16a}, \eqref{liu-16b}  (see the Appendix for the details)
obtained from the leading terms of \eqref{biham}. Such a  bihamiltonian structure is said to be of hydrodynamic type. 

It was shown by Dubrovin that to any Frobenius manifold $M$ one can associate an
integrable hierarchy, called the {\em principal hierarchy} of $M$, consisting of evolutionary PDEs of the form \eqref{principal} given in the Appendix, which
possesses a bihamiltonian structure \eqref{liu-16a}, \eqref{liu-16b}  of hydrodynamic type \cite{Du2}. For any semisimple Frobenius manifold, Dubrovin and Zhang provide an axiomatic approach to 
define a unique deformation of the principal hierarchy of the form \eqref{liu-12}, called the topological deformation. The bihamiltonian structure \eqref{liu-16a}, \eqref{liu-16b} also has a deformation \eqref{liu-18a}, \eqref{liu-18b} which is conjectured to have the polynomial property possessed by the bihamiltonian structure \eqref{biham} of the Drinfeld-Sokolov hierarchy. Note that the polynomial property of the first Poisson bracket \eqref{liu-18a} is proved in \cite{BPS1, BPS2}. 

For a given simple Lie algebra $\g$ of type $X_n$ with an appropriate choice of a basis of the principal Heisenberg subalgebra of $\hg$, it follows from Theorem 4.3 of \cite{DLZ2} that the dispersionless Drinfeld-Sokolov hierarchy together with its bihamiltonian structure of hydrodynamic type coincides 
with the principal hierarchy of the Frobenius manifold $M$ defined on the orbit space
of the Coxeter group of type $X_n$. When the Lie algebra $\g$ is of ADE type, this Frobenius manifold structure coincides with the one defined on the space of miniversal deformations of the ADE singularity \cite{S1}.

To see the relation between the Drinfeld-Sokolov hierarchy and the topological deformation of the principal hierarchy,  we introduce the normal coordinates $w_1, \dots, w_n$ of the Drinfeld-Sokolov hierarchy by using the densities of the Hamiltonians as follows:
\begin{equation}\label{norm-cor}
w_k=h_{j_k}(u,u_x,\dots),\quad k=1,\dots,n. 
\end{equation}
Here $j_1, \dots, j_n\in E_+$ are the first $n$ positive exponents in $E_+$, which coincide with
the exponents $m_1, \dots, m_n$ of $\g$ (see \S 14.2, \S 14.3 of \cite{Kac}).
This notion of normal coordinates is introduced in \cite{DZ3} for bihamiltonian integrable hierarchies with tau-structures. In terms of the these coordinates, the Drinfeld-Sokolov hierarchy is represented in the form \eqref{liu-12} after the introduction of the deformation parameter $\hbar$ as in \eqref{semi-classical}. Then by using the result of \cite{Wu2}
that the Virasoro symmetries of an ADE Drinfeld-Sokolov hierarchy act linearly on the tau functions associated to the tau structure $\Omega_{i, j}$ of the hierarchy, we know that the
ADE Drinfeld-Sokolov hierarchy coincides with the topological deformation of the principal hierarchy of the Frobenius manifold associated to the ADE singularity.  See the Appendix for a more detailed explanation.

For a simple Lie algebra of BCFG type,  the associated Drinfeld-Sokolov hierarchy, however, is not equivalent to the topological deformation of the principal hierarchy of the Frobenius manifold corresponding to the Coxeter group of BCFG type. This fact follows from the difference of the central invariants of the bihamiltonian structures of these two integrable hierarchies. The notion of central invariants of a semisimple bihamiltonian structure of the form  \eqref{biham} was introduced in \cite{DLZ1, LZ1}. Two bihamiltonian structures with the same leading terms are equivalent under Miura type transformations (see \cite{DZ3, DLZ1}) if and only if they have the same central invariants, which are $n$ functions of one variable. Here is a list of the central invariants of the Drinfeld-Sokolov bihamiltonian structures: 
\begin{equation}\label{ci-table}
\begin{array}{lccccc}
\g & & c_1 & \dots & c_{n-1} & c_n \\
 & & & & & \\
 & & & & & \\
A_n & & \frac1{24} & \dots & \frac1{24} & \frac1{24} \\
 & & & & & \\
 B_n &  &\frac1{24} & \dots & \frac1{24} & \frac1{12} \\
 & & & & & \\
 C_n & & \frac1{12} & \dots & \frac1{12} & \frac1{24} \\
 & &  & & & \\
 D_n & & \frac1{24} & \dots & \frac1{24} & \frac1{24} \\
 & & & & & \\
E_n, ~n=6, \, 7,\, 8 &  & \frac1{24} & \dots & \frac1{24} & \frac1{24} \\
 & & & & & \\
F_n, ~ n=4 & & \frac1{24} & \frac1{24} & \frac1{12} & \frac1{12} \\
 & & & & & \\
 G_n, ~ n=2 & & \frac1{8} & & & \frac1{24}
\end{array}
\end{equation}
Note that if the bilinear form $(\ |\ )$ that we introduced in the beginning of Section \ref{sec-4-1} is not the normalized one, the central invariants have the following formula:
\[c_i=\frac{(\alpha_i^\vee\, |\, \alpha_i^\vee)}{48},\quad i=1,\dots,n,\]
where $\alpha_i^\vee$ are simple coroots of $\g$. We also note that the central invariants are ordered according to the order of the canonical coordinates
of the Frobenius manifold. The results in the table \eqref{ci-table} correspond
to a specific ordering of the canonical coordinates. 

On the other hand, it is shown on page 832 of \cite{DLZ2} that if a bihamiltonian structure is equivalent under a Miura type transformation to the bihamiltonian structure of the topological 
deformation of the principal hierarchy of a semisimple Frobenius manifold, then all its central invariants must be equal to $\frac1{24}$. Thus, from the list of the central invariants of the BCFG Drinfeld-Sokolov bihamiltonian structures given in \eqref{ci-table}, we conclude that the BCFG Drinfeld-Sokolov hierarchies are not equivalent to the topological deformation of the principal hierarchy of any semisimple Frobenius manifold.
Nonequivalent deformations cannot possess a common solution, because their genus one free energies
do not coincide (see the next subsection for examples). Thus Theorem \ref{jw-3} is proved.

\subsection{The $\Gamma$-reduction theorem}\label{liu-22}

Given a simple Lie algebra $\g$ together with a group $\Gamma$
generated by an automophism $\sigma$ of $\g$, we show in this subsection that one can carefully choose a gauge in the Drinfeld-Sokolov reduction procedure to ensure the existence of a $\sigma$-action
on the resulting Drinfeld-Sokolov hierarchy. If the subalgebra  $\g^{\sigma}\subset\g$ given by the fixed points of $\sigma$ is also a simple Lie algebra, then we show that the flows of the Drinfeld-Sokolov hierarchy that are fixed by the action
of $\sigma$ yield exactly the Drinfeld-Sokolov hierarchy corresponding to $\g^\sigma$. Furthermore, by restricting the tau functions of the Drinfeld-Sokolov hierarchy associated to $\g$ to the fixed points of the action of $\sigma$ on the big phase space, we provide a way to obtain tau functions of the Drinfeld-Sokolov hierarchy associated to $\g^\sigma$.

Let us specify $\g$ to be a simple Lie algebra with one of the Dynkin diagrams given in Section \ref{introd}.
It is well-known that the outer automorphisms of $\g$ are induced by automorphisms of its Dynkin diagram, which we listed in \eqref{jw-12-21}--\eqref{jw-12-24}.
Suppose the rank of $\g$ is $n$ and $\{E_i, H_i, F_i\}_{i=1}^n$ is a set of Weyl generators of $\g$.
Then the induced automorphism of $\g$ is obtained from the following relations:
\begin{equation}\label{jw-12-25}
\sigma(E_i)=E_{\bar{\sigma}(i)}, \quad \sigma(F_i)=F_{\bar{\sigma}(i)}, \quad i=1, \dots, n.
\end{equation}
The automorphism $\sigma$ of $\g$ can be extended to an automorphism of $\hg$ which acts on the loop variable
$\lambda$ and on the central element $K$ trivially.

\begin{pro}[\S 7.9 of \cite{Kac}]
Let $\g$ and $\sigma$ be given as above.  Then the set  $\g^{\sigma}\subset\g$ of fixed points of $\sigma$
is a simple Lie algebra of type $B_n$, $C_n$, $F_4$, $G_2$ respectively.
\end{pro}
Kac's proof is based on an explicit description of a Chevalley basis of $\g$ which has nice properties under the action of
$\sigma$. By using his method, one can also show that $\sigma(E_0)=E_0$ and $\sigma(F_0)=F_0$. In particular
$\sigma(\Lambda)=\Lambda$, where $\Lambda$ is defined in \eqref{zh-Lambda}. So the principal Heisenberg subalgebra $\s=\Ker\ad_\Lambda$ is an invariant
subspace of $\sigma$.

It is easy to verify that in all four cases, $\sigma(I_{\pm})=I_\pm$ and $\sigma(\rho)=\rho$, so the subspaces $\Ker\ad_{I_\pm}$ are invariant under the action of $\sigma$. Note that the homogeneous
component $\g\cap\hg^j$ is the eigen-subspace of $\ad_\rho$ with eigenvalue $j$ and is also an invariant subspace of $\sigma$.

\begin{lem}\label{lem-54}
One can choose the lowest weight vectors $\gamma_1, \dots, \gamma_n$ to be eigenvectors of $\sigma$.
\end{lem}
\begin{prf}
Note that $\gamma_i\in \Ker \ad_{I_-} \cap\,\g\cap \hg^{-m_i}$, and the two spaces $\Ker \ad_{I_-}$ and
$\g\cap\hg^{-m_i}$ are both invariant under the action of $\sigma$. We can take $\gamma_i$ to be the
eigenvectors of $\sigma$.

In fact, except the $D_{2\ell}$ case, all the intersections $\Ker \ad_{I_-} \cap\,\g\cap \hg^{-m_i}$ have dimension one,
and the $\gamma_i$ are eigenvectors automatically. In the $D_{2\ell}$ case, the exponents $m_{\ell}=m_{\ell+1}=2\ell-1$.
We need to carefully choose $\gamma_{\ell}$ and $\gamma_{\ell+1}$ such that they are eigenvectors of $\sigma$.
\end{prf}

Suppose $\sigma(\gamma_i)=z_i \gamma_i$. Note that the order of $\sigma$ is $2$ (in Case 1, 2, 3) or $3$ (in Case 4),
and $z_i$ can be $\pm 1$ (in Case 1, 2, 3) or $1$, $\omega$, $\omega^2$ (in Case 4), where $\omega=e^{2\,\pi\,i/3}$.

According to  definition \eqref{eq-Li}, the action of $\sigma$ on the irreducible module $L_i$ is the same as the action of 
$z_i\Id$, so we have
\[\g^\sigma=\bigoplus_{\{i\mid z_i=1\}}L_i.\]

\begin{lem}
The subspace $\mathrm{Span}\{\gamma_i\mid z_i=1\}$ gives the lowest weight gauge of $\g^\sigma$.
\end{lem}
\begin{prf}
The Weyl generators $\{E_i^\sigma, H_i^\sigma, F_i^\sigma\}_{i=1}^n$ of $\g^\sigma$ can be chosen as
(see \S 8.3 of \cite{Kac})
\begin{align*}
D_{n+1}: \ &K_i^\sigma=K_i\ (1\le i \le n-1),\quad K_n^\sigma=K_n+K_{n+1}; \\
A_{2n-1}: \ &K_i^\sigma=K_i+K_{2n-i}\ (1 \le i \le n-1), \quad K_n^\sigma=K_n; \\
E_6: \ &K_i^\sigma=K_{2i-1}+K_{7-i}\ (i=1, 2),\quad K_j^\sigma=K_{2j-4}\ (j=3, 4); \\ 
D_4: \ &K_1^\sigma=K_1+K_3+K_4,\quad K_2^\sigma=K_2,
\end{align*}
where $K=E, H, F$. It is easy to see that the generators $I_+^\sigma, 2\rho^\sigma, I_-^\sigma$ of the $sl_2$ subalgebra of $\g^\sigma$ coincide respectively with $I_+, 2\rho, I_-$.
Hence we have,
for $i$ such that $z_i=1$,
\[[\rho^\sigma, \gamma_i]=-m_i\gamma_i, \quad [I_-^\sigma, \gamma_i]=0.\]
The lemma is proved.
\end{prf}

Let $\Lcal^{\can}(u^i)$ be the  differential operator chosen according to the lowest weight gauge.
We have $\sigma(\Lcal^{\can}(u^i))=\Lcal^{\can}(z_iu^i)$. If we regard $\sigma: V\to V$ as
a smooth map from the manifold $V$ to $V$, then the above relation shows that the pullback of $\sigma$
acts on the coordinates $(u^1, \dots, u^n)$ by $\sigma^*(u^i)=z_i\, u^i$. Here $V$ is defined in \eqref{fix-gauge}. We extend this action
to the jet space of $V$ by $\sigma^*(u^{i,s})=z_i\,u^{i,s}$. 
Next we need to study the following questions: (i) if $f$ is a differential polynomial of
$(u^1, \dots, u^n)$, what is $\sigma^*(f)$? (ii) if $X$ is a vector field on $J^\infty(V)$, what is
$d\sigma(X)$?

\begin{lem}
One can choose the basis $\{\Lambda_j\}_{j\in E}$ to be eigenvectors of $\sigma$.
\end{lem}
\begin{prf}
Similar to the proof of Lemma \ref{lem-54}.
\end{prf}

\begin{pro}\label{pro-57}
Suppose $\sigma(\Lambda_j)=\zeta_j\,\Lambda_j$.  Then for $j,k\in E_+$, we have 
\begin{itemize}
\item[i)] $\sigma^*(h_j)=\zeta_j^{-1} h_j$;
\item[ii)] $d\sigma(\frac{\p}{\p t^j})=\zeta_j \frac{\p}{\p t_j}$;
\item[iii)] $\sigma^*(\Omega_{jk})=\zeta_j^{-1}\zeta_k^{-1}\Omega_{jk}$;
\item[iv)] If $\Fcal$ is determined by \eqref{zh-12-1} and the string equation, then $\sigma^*\Fcal=\Fcal$.
\end{itemize}
\end{pro}
\begin{prf}
Let $U$ and $H$ be the two matrices appearing in Proposition \ref{pro-51} for $\Lcal^{\can}$.
Denote $\Lcal^{\can}$ by $\Lcal$ for short.
By using the fact that $\sigma(\Lcal(u))=\Lcal(\sigma^*(u))$ and the uniqueness of $U$ and $H$,
we obtain $\sigma(U)=U(\sigma^*(u))$ and $\sigma(H)=H(\sigma^*(u))$.  The fact that $H=-\sum_{j\in E_+}h_j\,\Lambda_{-j}/j$ implies that $\sigma(H)=-\sum_{j\in E_+}h_j\,\zeta_{-j}\,\Lambda_{-j}/j$.
Using the fact that $\zeta_{-j}=\zeta_j^{-1}$, this proves the first property of the proposition.

To prove property ii), let us consider the action of $\sigma$ on the following equation:
\[\frac{\p \Lcal}{\p t_j}=-[\left(e^{\ad_U}(\Lambda_j)\right)_+, \Lcal].\]
The action on the left hand side is
\[\sigma\left(\frac{\p \Lcal}{\p t_j}\right)=\sum_{i=1}^N\frac{\p \sigma^*(u^i)}{\p t_j}\gamma_i,\]
while the action on the right-hand side yields
\begin{align*}
&\sigma\left(-[\left(e^{\ad_U}(\Lambda_j)\right)_+, \Lcal]\right)
=-[\left(e^{\ad_{\sigma(U)}}(\sigma(\Lambda_j))\right)_+, \sigma(\Lcal)]\\
=&-\zeta_j[\left(e^{\ad_{U(\sigma^*(u))}}(\Lambda_j)\right)_+, \Lcal(\sigma^*(u))]\\
=& \zeta_j \sigma^*\left(\sum_{i=1}^N\frac{\p u^i}{\p t_j}\gamma_i\right).
\end{align*}
Comparing the above two formulae, we prove the second property of the proposition.

The properties iii) and iv) are easy consequences of i) and ii). The proposition is proved.
\end{prf}

Now we  prove the $\Gamma$-reduction theorem.

\begin{thm}\label{liu-5}
(a) Each flow $\frac{\p}{\p t_j}$ of the Drinfeld-Sokolov hierarchy associated to the simple Lie algebra $\g$ with $\zeta_j=1$ can be restricted to the jet space
of $V^\sigma=\{q\in V|\sigma(q)=q\}$. Moreover, all these restricted flows yield the Drinfeld-Sokolov hierarchy associated to the simple Lie algebra
$\g^\sigma$.
(b) Let $\tau$ be a tau function of the Drinfeld-Sokolov hierarchy associated to $\g$, which is a function on the big phase space
$H_{\g}=\mathrm{Spec}(\CC[t_j\mid j\in E_+])$. Let $H_{\g}^\sigma=\mathrm{Spec}(\CC[t_j\mid j\in E_+,\ \zeta_j=1])$.
The restriction of $\tau$ to $H_{\g}^\sigma$ gives a tau function of the Drinfeld-Sokolov hierarchy associated to $\g^\sigma$.
\end{thm}
\begin{prf}
Proposition \ref{pro-57} implies that if $\zeta_j=1$, then $\frac{\p}{\p t_j}$ is invariant under the action of
$\sigma$. Hence it can be restricted to the set of invariant subspace $V^\sigma$.
Let $\s^\sigma$ be the principal Heisenberg subalgebra of $\hg^\sigma$.  $\{K\}\cup\{\Lambda_j\mid j\in E, \zeta_j=1\}$
defines a basis of $\s^\sigma$. Then all the flows of the Drinfeld-Sokolov hierarchy associated to $\g^\sigma$ can be obtained from that of
$\g$ and we  prove the first part of the theorem. The second part of the theorem is an easy consequence of the first part. 
\end{prf}

It is well-known that the genus one free energy of the ADE singularity is given by
\[\Fcal^1=\frac{1}{24}\log\det M,\]
where $M=(M^\alpha_\beta)$ with $M^\alpha_\beta=\eta^{\al\gamma}\frac{\p^3 \Fcal^0}{\p t^{1,0}\p t^{\gamma,0}\p t^{\beta,0}}$, and $(\eta^{\alpha\beta})=(\eta_{\al\beta})^{-1}$ with $\eta_{\al\beta}=\left.\frac{\p^3 \Fcal^0}{\p t^{1,0}\p t^{\al,0}\p t^{\beta,0}}\right|_{t^{\al,p\ge 1}=0}$.
It can also be written as
\[\Fcal^1=\frac{1}{24}\sum_{i=1}^n\log \frac{\p u_i(v(t))}{\p t^{1,0}}=\frac{1}{24}\sum_{i=1}^n\log u_{i,x},\]
where $u_i(v)$ is the $i$-th canonical coordinates (different from the functions $u^i(x)$ that appear in $\Lcal$), and $v(t)=(v^1(t),\dots, v^n(t))$ is the topological solution  of the principal hierarchy, described in the Appendix. For the relation of the variables $t^{\al,p}$ with the time variables $t_j, j\in E_+$ of the Drinfeld-Sokolov hierarchies and the
definition of $\Fcal^0$, see Section \ref{jw-12-1} and the Appendix. 

When restricting $\Fcal^1$ to the subspace $H^\sigma_\g$ defined in Theorem \ref{liu-5},
one can show that some canonical coordinates will become equal. The genus
one free energies of the $\sigma$-invariant sectors are given by
\begin{align*}
&D_{n+1}\rightsquigarrow B_n:\ \Fcal^1=\sum_{i=1}^{n-1}\frac{1}{24}\log u_{i,x}+\frac{1}{12}\log u_{n,x},\\
&A_{2n-1}\rightsquigarrow C_n:\ \Fcal^1=\sum_{i=1}^{n-1}\frac{1}{12}\log u_{i,x}+\frac{1}{24}\log u_{n,x},\\
&E_{6}\rightsquigarrow F_4:\ \Fcal^1=\sum_{i=1}^{2}\frac{1}{24}\log u_{i,x}+\sum_{i=3}^4 \frac{1}{12}\log u_{i,x},\\
&D_4\rightsquigarrow G_2:\ \Fcal^1=\frac{1}{8}\log u_{1,x}+\frac{1}{24}\log u_{2,x}.
\end{align*}
Note that the coefficients before $\log$'s are exactly the central invariants of the corresponding Drinfeld-Sokolov bihamiltonian structures.

These expressions manifest the fact that the total descendent potential of 
a semisimple cohomological field theory never satisfies
the BCFG Drinfeld-Sokolov hierarchies, since its genus one part has the expression
\[\Fcal^1=\sum_{i=1}^{n}\frac{1}{24}\log u_{i,x}+G(u),\]
where $G(u)$ is the $G$-function of the corresponding Frobenius manifold \cite{DZ1}.

\subsection{Proof of the main theorem}

We follow the idea spelled out in the beginning of the section.  After the $\Gamma$-reduction theorem, the only remaining part of the proof is to compare the $\Gamma$-action on genus zero part of FJRW-theory with its counterpart in integrable hierarchies.
We can use the mirror theorems of \cite{FJR2, FJMR} to compare the action of FJRW theory with the action of the mirror Frobenius manifold structures. In this
subsection, we match the $\Gamma$-action on the Frobenius manifold structure of the ADE singularities with that of the semiclassical limit of the Drinfeld-Sokolov hierarchy. The latter object is well studied in \cite{DLZ2}.

\vskip 1em

\noindent\textbf{$A_{2n-1}$ case:}

The polynomial $W$ is given by $z^{2n}$, and the action of $\sigma$ is $\sigma(z)=-z$. We redenote this polynomial by $f(z)$. The miniversal deformation $f_t(z)$ is chosen as
\begin{equation}
f_t(z)=z^{2n}+t_{2n-1}z^{2n-2}+\cdots+t_2\,z+t_1. \label{superpotential}
\end{equation} 
To ensure that $f_t$ is $\sigma$-invariant, we need to define
\begin{equation}
\sigma^*(t_i)=(-1)^{i-1}t_i, \label{sigma-t-A}
\end{equation}
which gives the action of $\sigma$ on the corresponding Frobenius manifold. It is easy to check that the mirror map of \cite{FJR2}
matches it with the corresponding action on FJRW-theory. 

The Drinfeld-Sokolov hierarchy in this case is controlled by a scalar operator
\[L=D^{2n}+\tilde{u}^{2n-1}\,D^{2n-2}+\cdots+\tilde{u}^2\,D+\tilde{u}^1.\]
The action of $\sigma$ is given by
\begin{equation}
\sigma(L)=L^\dagger=(-D)^{2n}+(-D)^{2n-2}\,\tilde{u}^{2n-1}+\cdots+D\,\tilde{u}^2+\tilde{u}^1. \label{sigma-L-A}
\end{equation}

By taking the semiclassical limit, the operator $L$ becomes the miniversal deformation \eqref{superpotential},
and every $\tilde{u}_i$ tends to $t_i$. Then the action \eqref{sigma-L-A} implies the relation \eqref{sigma-t-A}.

\vskip 1em

\noindent\textbf{$D_{n+1}$ case:}

The polynomial $W$ is given by $f(x,y)=x^{n-1}+xy^2$, and $\sigma^*(x,y)=(x,-y)$. The miniversal deformation has the expression
\[f_t(x,y)=P(x)+xy^2+t_n y,\quad \mbox{where } P(x)=x^{n-1}+t_{n-1}x^{n-2}+\dots+t_1.\]

\begin{lem}
Saito's Frobenius manifold structure for $f(x,y)$ can be obtained from the super potential
\begin{equation}
\lambda(z)=P(z^2)-\frac{t_n^2}{4z^2}\label{superpotential-D}
\end{equation}
and the residue pairing
\[\langle\p_1,\p_2\rangle=-2(\Res_{z=\infty}+\Res_{z=0})\frac{\p_1\lambda(z)\,\p_2\lambda(z)}{\lambda'(z)}.\]
\end{lem}
\begin{prf}
Denote by $H$ the Hessian determinant of $f_t$.  Then the original residue pairing can be defined as
\[\langle\p_1,\p_2\rangle=-\Res_{\{\infty\}} \frac{\p_1 f_t\,\p_2 f_t}{\p_x f_t\, \p_y f_t}dx\wedge dy=
\sum_{(x_i,y_i)}\left.\frac{\p_1f_t\,\p_2f_t}{H}\right|_{(x,y)=(x_i,y_i)},\]
where $(x_i, y_i)$ runs over all the critical points of $f_t$.

It is easy to see that the critical points of $f_t$ are given by
\[x_i=z_i^2,\quad y_i=-\frac{t_n}{2z_i^2},\]
where $z_i^2$ runs over all the critical points of $\tilde{\lambda}(x)=P(x)-\frac{t_n^2}{4x}$.
Note that at these critical points we have
\begin{align*}
&\left.f_t\right|_{(x,y)=(x_i,y_i)}=\left.\lambda(z)\right|_{z=z_i},\\
&\left.H\right|_{(x,y)=(x_i,y_i)}=\frac12\left.\lambda''(z)\right|_{z=z_i}.
\end{align*}
Then the lemma is proved by the residue theorem.
\end{prf}

The action of $\sigma$ is given by
\begin{equation}
\sigma^*(t_i)=t_i\ (1\le i \le n-1),\quad \sigma^*(t_n)=-t_n. \label{sigma-t-D}
\end{equation}
The super potential \eqref{superpotential-D} is trivially invariant under the action of this $\sigma$.
It is easy to check that the mirror map of \cite{FJR2}
matches it with the corresponding action on FJRW-theory.

The $D_{n+1}$ Drinfeld-Sokolov hierarchy is governed by the following scalar pseudo-differential operator (see \cite{DS, LWZ}):
\[L=D^{-1}\left(D^{2n+1}+\sum_{i=1}^n\left(\tilde{u}^i D^{2i-1}+D^{2i-1}\tilde{u}^i\right)+\rho D^{-1}\rho\right).\]
The action of $\sigma$ on the hierarchy is induced by
\[\sigma^*(\tilde{u}^i)=\tilde{u}^i,\quad \sigma^*(\rho)=-\rho,\]
whose semiclassical limit is exactly \eqref{sigma-t-D}.

The $D_4$ case with a $\ZZ_3$-action will be treated with full details in the next section, so we omit it here.

\vskip 1em

\noindent\textbf{$E_6$ case:}

The polynomial $W$ in this case is given by $f(x,y)=x^3+y^4$, and the miniversal deformation can be chosen as
\[f_t(x,y)=t_1+t_2\,y+t_3\,x+t_4\,y^2+t_5\,x\,y+t_6\,x\,y^2+x^3+y^4.\]
The action of $\sigma$ on $(x,y)$ is given by $(x,-y)$. We have
\begin{equation}
\sigma^*(t_i)=t_i\ (i=1,3,4,6),\quad \sigma^*(t_i)=-t_i\ (i=2,5). \label{sigma-E-t}
\end{equation}

There is no simple description of the Drinfeld-Sokolov hierarchy of $E$ type. We must start from
a simple Lie algebra of $E_6$ type, and perform the Drinfeld-Sokolov reduction from scratch. Note that a similar computation
has been done in \cite{DLZ2}, but the action of $\sigma$ was not considered in that paper.

The $E_6$ Lie algebra has a matrix realization of dimension $27$,
which can be generated by the following Weyl generators:
\begin{align*}
E_1&=e_{6,7}+e_{8,9}+e_{10,11}+e_{12,14}+e_{15,17}+e_{26,27},\\
E_2&=e_{4,5}+e_{6,8}+e_{7,9}-e_{18,20}-e_{21,22}-e_{23,24},\\
E_3&=e_{4,6}+e_{5,8}+e_{11,13}+e_{14,16}+e_{17,19}+e_{25,26},\\
E_4&=e_{3,4}-e_{8,10}-e_{9,11}-e_{16,18}-e_{19,21}+e_{24,25},\\
E_5&=e_{2,3}-e_{10,12}-e_{11,14}-e_{13,16}+e_{21,23}+e_{22,24},\\
E_6&=e_{1,2}+e_{12,15}+e_{14,17}+e_{16,19}+e_{18,21}+e_{20,22},\\
F_i&=E_i^t,\quad H_i=[E_i,F_i],
\end{align*}
where $e_{i,j}$ denotes the matrix with $(i,j)$-th entry being $1$ and other entries vanishing, and $A^t$ means the transpose of $A$.
The action of $\sigma$ is induced by the action on generators
\begin{align*}
&\sigma(K_1)=K_6,\quad \sigma(K_3)=K_5,\quad \sigma(K_2)=K_2,\\
&\sigma(K_6)=K_1,\quad \sigma(K_5)=K_3,\quad \sigma(K_4)=K_4,
\end{align*}
where $K=E, F, H$. 

To perform the Drinfeld-Sokolov reduction, we need to take a principal nilpotent $sl_2$ subalgebra and consider the decomposition
of $\g$ as an $sl_2$-module. The $sl_2$ subalgebra is given by $\{I_+, 2\rho, I_-\}$, where (according to Section \ref{sec-411}):
\begin{align*}
I_+=&E_1+E_2+E_3+E_4+E_5+E_6,\\
\rho=&8 H_1+11 H_2+15 H_3+21 H_4+15 H_5+8 H_6,\\
I_-=&16 F_1+22 F_2+30 F_3+42 F_4+30 F_5+16 F_6.
\end{align*}
The lowest weight vectors of the $sl_2$-module $\g$ can be chosen as
\begin{align*}
\gamma_1=&\left(F_1+F_6\right)+\frac{15}{8}\left(F_3+F_5\right)+\frac{11}{8}F_2+\frac{21}{8}F_4,\\
\gamma_2=&F_{111100}-F_{010111}+\frac{15}{11}\left(F_{101110}+F_{001111}\right),\\
\gamma_3=&F_{111110}-F_{011111}+\frac{16}{11}F_{101111}+\frac{21}{8}F_{011210},\\
\gamma_4=&F_{112210}+F_{011221}+\frac{8}{15}F_{111211},\\
\gamma_5=&F_{112211}+F_{111221},\\
\gamma_6=&F_{122321},
\end{align*}
where $F_{n_1n_2n_3n_4n_5n_6}$ is the basis of the root space $\g_{-\alpha}$ with
\[\alpha=n_1\alpha_1+n_2\alpha_2+n_3\alpha_3+n_4\alpha_4+n_5\alpha_5+n_6\alpha_6,\]
such that the first nonzero entry of the first nonzero column of $F_{n_1n_2n_3n_4n_5n_6}$ is equal to $1$.

The action of $\sigma$ on these lowest weight vectors is given by
\[\sigma(\gamma_i)=\gamma_i\ (i=1,3,4,6),\quad \sigma(\gamma_i)=-\gamma_i\ (i=2,5).\]
Let $u^i\ (i=1, \dots,6)$ be the coordinates of $V$ w.r.t. the basis $\gamma_i\ (i=1, \dots, 6)$, then we have
\begin{equation}
\sigma^*(u^i)=u^i\ (i=1,3,4,6),\quad \sigma^*(u^i)=-u^i\ (i=2,5). \label{sigma-E-u}
\end{equation}

The last step is to relate $t_i$ and $u^i$. Note that they are different coordinate systems on the same Frobenius
manifold. Hence, we can compare them with a fixed collection of flat coordinates. Let $v^i\ (i=1, \dots, 6)$ be a system of flat coordinates.
The coordinates $t_i, u_i, v^i$ have the following degrees:
\begin{align*}
&\deg t_1=\deg u^6=\deg v^1=12,\quad \deg t_2=\deg u^5=\deg v^2=9,\\
&\deg t_3=\deg u^4=\deg v^3=8, \quad \deg t_4=\deg u^3=\deg v^4=6, \\
&\deg t_5=\deg u^2=\deg v^5=5, \quad \deg t_6=\deg u^1=\deg v^6=2.
\end{align*}
Hence $v^i$ must be polynomials of $t_i$ or $u^i$, since they are homogeneous holomorphic
functions with positive degrees.
Note that $\sigma$ only changes the sign of $t_i$ or $u^i$ with odd degree. Namely, if $v^i$ has even degree,
then it must be $\sigma$-invariant, while $\sigma^*(v^i)=-v^i$ for odd degree $v^i$.
Thus we have shown that both $\sigma$-actions on the Frobenius manifold that are induced from the $\sigma$-action \eqref{sigma-E-t} on $t$ and from the $\sigma$-action \eqref{sigma-E-u} on $u$ are the same.
It is easy to check that the mirror map of \cite{FJR2}
matches it with the corresponding action on FJRW-theory. Thus we proved Theorem \ref{liu-21}. 


\section{Examples}\label{jw-12-1}

The Drinfeld-Sokolov hierarchies together with the string equation determine the free energies $\Fcal^g$. In  practice, it takes
considerable work to use the integrable hierarchies to do explicit calculation. In this section, we give an algorithm to obtain the FJRW invariants from the corresponding 
integrable systems for three examples. 
We first consider the $D_4$ Drinfeld-Sokolov hierarchy and its free energies. Then we consider the $B_3$ and $G_2$ Drinfeld-Sokolov hierarchies,
these two integrable hierarchies and their free energies can be obtained from that of the $D_4$ case by taking the invariant loci of
certain symmetric group actions. Note that the algorithm to compute the free energies, which is designed for the $D_4$ case, can also be applied directly to $B_3$ and $G_2$ Drinfeld-Sokolov hierarchies. Since the reduced integrable hierarchies are simpler then that of the $D_4$ case, the direct computation is more efficient.

\subsection{The $D_4$ Drinfeld-Sokolov hierarchy} \label{sec-51}

The simple Lie algebra of type $D_4$ is just $\mathfrak{o}(8)$.
We choose, as in \cite{DS}, the following matrix realization of this Lie algebra: 
\[\mathfrak{o}(8)=\{A\in \mathrm{gl}(8,\mathbb{C})| A=-S A^T S^{-1}\},\]
where $A^T$ denotes the transposition with respect to the secondary diagonal,
and $S=\textrm{diag}(1,-1,1,-1,-1,1,-1,1)$.

The Weyl generators of $\g$ can be chosen as
\begin{align*}
&E_i=e_{i+1,i}+e_{9-i,8-i}, \quad E_4=\frac12 (e_{5,3}+e_{6,4}),\\
&F_i=e_{i,i+1}+e_{8-i,9-i}, \quad F_4=2 (e_{3,5}+e_{4,6}),\\
&H_i=-e_{i,i}+e_{i+1,i+1}-e_{8-i,8-i}+e_{9-i,9-i},\\
&H_4=-e_{3,3}-e_{4,4}+e_{5,5}+e_{6,6},
\end{align*}
where $1\le i\le 3$, and $e_{i,j}$ is the matrix with its $(i,j)$-th entry being $1$ and other entries vanishing.
We take a set of Chevalley bases
\[K_5:=K_{1,2},\ K_6:=K_{3,2},\ K_7:=K_{4,2},\]
\[K_8:=K_{5,3}=K_{6,1},\ K_9:=K_{6,4}=K_{7,3},\ K_{10}:=K_{7,1}=K_{5,4},\]
\[K_{11}:=K_{1,9}=K_{3,10}=K_{4,8},\ K_{12}:=K_{11,2}=K_{5,9}=K_{6,10}=K_{7,8},\]
where $K_{i,j}=[K_i,K_j]$, and $K$ denotes either $E$ or $F$.

The $sl_2$ subalgebra $\{I_+, 2\rho, I_-\}$ reads
\begin{align*}
I_+&=E_1+E_2+E_3+E_4,\\
2\rho&=6H_1+10H_2+6H_3+6H_4,\\
I_-&=6F_1+10F_2+6F_3+6F_4.
\end{align*}
The lowest weight vectors $\gamma_1, \dots, \gamma_4$ can be chosen as
\begin{align*}
&\gamma_1=3F_1+5F_2+3F_3+3F_4,\ \gamma_4=F_{12},\\
&\gamma_2=F_8+F_{10}-2F_9,\ \gamma_3=F_8-F_{10}.
\end{align*}
We denote the automorphisms \eqref{jw-12-25} of $\g$  induced by \eqref{jw-12-21} and
\eqref{jw-12-24} by $\sigma_1$ and $\sigma_4$ respectively.
One can check that
\begin{equation}
\sigma_1(\gamma_1)=\gamma_1,\ \sigma_1(\gamma_2)=\gamma_2,\ \sigma_1(\gamma_3)=-\gamma_3,\ \sigma_1(\gamma_4)=\gamma_4.\label{sigma-1}
\end{equation}
These lowest weight vectors are eigenvectors of $\sigma_1$, and  will be used in our reduction procedure from the $D_4$  Drinfeld-Sokolov hierarchy to the $B_3$ case.
For the reduction of the $D_4$ Drinfeld-Sokolov hierarchy to the $G_2$ case, we choose the following base of the space of  the lowest weight vectors: 
\begin{equation}\label{lf-1}
\gamma'_1=\gamma_1,\ \gamma'_2=F_8+\omega F_9+\omega^2 F_{10},\ \gamma'_3=F_8+\omega^2 F_9+\omega F_{10},\ \gamma'_4=\gamma_4,
\end{equation}
where $\omega=\frac{-1+\sqrt{-3}}{2}$. They are the eigenvectors of $\sigma_4$, as the following relations show:
\begin{equation}
\sigma_4(\gamma'_1)=\gamma'_1,\ \sigma_4(\gamma'_2)=\omega^2\gamma'_2,\ \sigma_4(\gamma'_3)=\omega\gamma'_3,\ \sigma_4(\gamma'_4)=\gamma'_4.\label{sigma-4}
\end{equation}
In what follows, we will illustrate in detail the reduction procedure from the $D_4$ case to the $B_3$ case by using the lowest weight vectors $\gamma_1, \dots, \gamma_4$. The reduction procedure from $D_4$ case to the $G_2$ case is similar.

The matrix differential operator $\Lcal^{\can}$ has the form
\begin{equation}
\Lcal^{\can}=\frac{d}{dx}+\Lambda+\sum_{i=1}^4 u^i(x)\gamma_i,
\end{equation}
where $\Lambda=I_+-\lambda \gamma_4$.
The associated Drinfeld-Sokolov hierarchy can be represented, as shown  in \cite{DS}, by using the following pseudo-differential operator:
\begin{equation}\label{lax-scalar}
L=\p_x^6+\p_x^{-1}\sum_{i=1}^3(\tilde{u}^i \p_x^{2i-1}+\p_x^{2i-1} \tilde{u}^i)+\p_x^{-1}\tilde{u}^4\,\p_x^{-1}\tilde{u}^4,
\end{equation}
where
\begin{align}
\tilde{u}^1=& 2u^4+28u^1u^2-144(u^1)^3-\frac{293}{2}(u^1_x)^2 \nonumber\\
& -118 u^1u^1_{xx}+4u^2_{xx}-17 u^1_{xxxx},\label{tildu-1}\\
\tilde{u}^2=& 98(u^1)^2+28u^1_{xx}-6u^2, \label{tildu-2}\\
\tilde{u}^3=& -14u^1,\quad \tilde{u}^4=2u^3.\label{tildu-3}
\end{align}

The dispersionless limit of \eqref{lax-scalar} gives the super potential of the corresponding Frobenius manifold.
We write this super potential in the following form:
\begin{align}
&\lambda(p)=p^6+ v^3 p^4+\left( v^2+\frac14 (v^3)^2\right) p^2+ v^1+\frac16 v^2 v^3\nonumber \\
&\qquad \quad +\frac1{108} (v^3)^3+\frac{(v^4)^2}{4 p^2}.\label{zh-11-16}
\end{align}
The structure constants of the Frobenius manifold can be obtained by
\[c_{ijk}=-\Res_{p=\infty} \frac{\p_i \lambda(p)\p_j \lambda(p)\p_k \lambda(p)}{\lambda'(p)}-\Res_{p=0} \frac{\p_i \lambda(p)\p_j \lambda(p)\p_k \lambda(p)}{\lambda'(p)},\]
where $\p_i=\frac{\p}{\p v^i}$, and $v^1,\dots, v^4$ are the flat coordinates. The potential of the Frobenius manifold reads
\begin{align}
F(v)=&\frac1{12}v^1(v^2)^2+\frac1{12}(v^1)^2v^3+\frac14v^1 (v^4)^2+\frac{1}{24}v^2 v^3 (v^4)^2
\nonumber\\
&-\frac{(v^2)^3v^3}{216}+\frac{(v^3)^3(v^4)^2}{432}+\frac{(v^2)^2(v^3)^3}{1296}
+\frac{(v^3)^7}{1632960}, \label{zh-11-3}
\end{align}
and the Euler vector field is given by
\begin{equation}
E=v^1 \frac{\p}{\p v^1}+\frac23 v^2\frac{\p}{\p v^2}+\frac13 v^3 \frac{\p}{\p v^3}+\frac23 v^4\frac{\p}{\p v^4}.\label{zh-11-17}
\end{equation}
They coincide with the one obtained from the FJRW theory of $D_4$ singularity \cite{FJMR} (see also \cite{LYZ} for another approach by using $W$-constraints). Note 
that there is a sign difference between the super potential $\lambda(p)$ given in 
\eqref{zh-11-16} and the one that is used in \cite{FJMR}. Due to this we need to
make the change $(v^4)^2\to -(v^4)^2$ in order to identify the above potential $F$ with the one given in \cite{FJMR}. We keep using the super potential \eqref{zh-11-16} for the convenience of the presentation of the integrable hierarchies.

To write down the full Drinfeld-Sokolov hierarchy associated to the scalar operator \eqref{lax-scalar},
we introduce two operators
\begin{align*}
P=&L^{1/6}=\p_x+\sum_{i\ge 1}P_i \p_x^{-i},\\
Q=&L^{1/2}=\p_x^{-1} \tilde{u}^4+\sum_{i\ge 0}Q_i \p_x^i.
\end{align*}
The operator $P$ is the sixth root of $L$ in the space of pseudo-differential operators of the first type,
while the operator $Q$ is the square root of $L$ in the space of pseudo-differential operators of the
second type defined in \cite{LWZ}. Note that both operators are convergent in their space. For more details, we refer the readers to
\cite{LWZ}.

The associated Drinfeld-Sokolov hierarchy is then given by
\begin{align}
\frac{\p L}{\p t_{k}}=[(P^k)_+, L], \quad k=1, 3, 5, \dots,\label{zh-11-1}\\
\frac{\p L}{\p \hat{t}_{k}}=[(Q^k)_+, L], \quad k=1, 3, 5, \dots.\label{zh-11-2}
\end{align}
The time variables $t^{\alpha,p}$ of the topological deformation have the expressions
\begin{equation}\label{zh-11-5}
t^{\alpha,p}=\left\{\begin{array}{cc}
\frac{6\Gamma\left(p+1+\frac{2\alpha-1}{6}\right)}{\Gamma\left(\frac{2\alpha-1}{6}\right)} t_{6p+2\alpha-1}, & \alpha=1, 2, 3; \ p\ge 0;\\ & \\
\frac{2 \Gamma\left(p+1+\frac12\right)}{\Gamma\left(\frac12\right)}\hat{t}_{\, 2p+1}, & \alpha=4;\  p\ge 0.
\end{array}\right.
\end{equation}
As defined in \cite{DZ3} and explained in \eqref{norm-cor}, the normal coordinates of this integrable hierarchy read
\[w^\al=\frac{6}{7-2\al}\res P^{7-2\al}\ (\al=1, 2, 3),\quad w^4=2\res Q.\]
They are the same as $\eta^{\alpha\beta}h_{\beta,-1}$, where
\begin{equation}\label{zh-11-6}
(\eta^{\alpha\beta})=\left(\begin{array}{cccc}
0 & 0 & 6 & 0 \\
0 & 6 & 0 & 0 \\
6 & 0 & 0 & 0 \\
0 & 0 & 0 & 2
\end{array}\right)
\end{equation}
is the contravariant metric of the corresponding Frobenius manifold, and $h_{\alpha,-1}$  are given in
\cite{LWZ} (page 1970).
More precisely,
\begin{align*}
w^1=&-\frac{2288}{27}(u^1)^3+4u^4+99 (u^1_x)^2+\frac{1012}9u^1u^1_{xx}-\frac{242}{45}u^1_{xxxx},\\
w^2=&-12u^2,\quad w^3=-28 u^1, \quad w^4=4u^3.
\end{align*}
To represent the pseudo-differential operator $L$ in terms of $w^1,\dots, w^4$, we need to use the
following relations:
\begin{align*}
&\tilde{u}^1=\frac12 w^1+\frac1{12} w^2 w^3+\frac1{216} (w^3)^3-\frac14 (w^3_x)^2\\
&\qquad -\frac13 w^2_{xx}-\frac29 w^3 w^3_{xx}+\frac{23}{45}w^3_{xxxx},\\
&\tilde{u}^2=\frac12 w^2+\frac18 (w^3)^2-w^3_{xx},\quad \tilde{u}^3=\frac12 w^3,\quad\tilde{u}^4=\frac12 w^4.
\end{align*}

Note that, in the above construction, the bilinear form $(\ \mid\ )$ is not the normalized one but differs by a factor $2$, and 
the central invariants $c_1=\dots=c_4=\frac1{12}$ \cite{LWZ}. To obtain the topological deformation, we need to rescale the
string coupling constant $\hbar$. More precisely, we first obtain the expressions of the flows 
\begin{align}
\frac{\p w^\al}{\p t^{\beta,q}},\quad \al, \beta=1,\dots,4,\quad  p, q\ge 0
\end{align}
from \eqref{zh-11-1}, \eqref{zh-11-2} and \eqref{zh-11-5}, and then perform the substitution 
\begin{equation}\label{zh-11-19}
\p_{t^{\beta,q}} w^\al\to \left(\frac{\hbar}{2}\right)^{\frac{1}2} \p_{t^{\beta,q}} w^\al,\quad
\p_x^k w^\al\to \left(\frac{\hbar}{2}\right)^{\frac{k}2} \p_x^k w^\al 
 \end{equation}
for $\al, \beta=1,\dots,4,\, q\ge 0, \, k\ge 1$.
The first few flows of the hierarchy read
\begin{align*}
\frac{\p w^\al}{\p t^{1,0}}=&w^\al_x,\quad \al=1,\dots,4;\\
\frac{\p w^1}{\p t^{2,0}}=&\frac{1}{18} w^2 w^3 w^3_x+\frac{1}{2}w^4 w^4_x
-\frac{1}{6} w^2 w^2_x+\frac{1}{36} (w^3)^2  w^2_x\\
&+\hbar \left(\frac{1}{12} w^3 w^2_{xxx}+\frac{1}{6}  w^2_{xx} w^3_x+\frac{5}{36} w^2_x w^3_{xx}+\frac{1}{18}  w^2 w^3_{xxx}\right)+\frac{\hbar^2}{20}  w^{2,5};\\
\frac{\p w^2}{\p t^{2,0}}=&w^1_x+\frac{1}{36} (w^3)^2 w^3_x-\frac{1}{6}w^2 w^3_x-\frac{1}{6} w^3 w^2_x\\
&+\hbar\left(\frac{1}{36} w^3_x w^3_{xx}+\frac{1}{36}  w^3 w^3_{xxx}
-\frac{1}{6}  w^2_{xxx}\right)+\frac{\hbar^2}{180}  w^{3,5};\\
\frac{\p w^3}{\p t^{2,0}}=&w^2_x;\\
\frac{\p w^4}{\p t^{2,0}}=&\frac{1}{6} w^3 w^4_x+\frac{1}{6} w^4 w^3_x+\frac{\hbar}{6}w^4_{xxx};
\end{align*}

\begin{align*}
\frac{\p w^1}{\p t^{3,0}}=&\frac{1}{12} (w^4)^2 w^3_x+\frac{1}{6} w^3 w^4 w^4_x
+\frac{1}{36} (w^2)^2 w^3_x+\frac1{1296} (w^3)^4 w^3_x+\frac{1}{18} w^2 w^3 w^2_x\\
&+\hbar\left( \frac{1}{108} w^3 (w^3_x)^3+\frac{1}{18} w^3 w^1_{xxx}+\frac{1}{12}  w^2_x w^2_{xx}
+\frac{1}{18} w^2 w^2_{xxx}+\frac{1}{12}  w^1_{xx} w^3_x\right.\\
&\quad \left.+\frac{1}{36} w^1_x w^3_{xx}
+\frac{1}{54}  (w^3)^2 w^3_x w^3_{xx}+\frac{1}{324} (w^3)^3 w^3_{xxx}
+\frac{1}{4} w^4_x w^4_{xx}+\frac{1}{6} w^4 w^4_{xxx}\right)\\
&+\hbar^2 \left(\frac{59}{1296} w^3_x (w^3_{xx})^2+\frac{2}{45} w^{1,5}
+\frac{7}{216} (w^3_x)^2 w^3_{xxx}+\frac{25}{648} w^3 w^3_{xx}w^3_{xxx}\right.\\
&\quad \left.+\frac{5}{216} w^3 w^3_x w^{3,4}+\frac{5}{1296} (w^3)^2 w^{3,5} \right)\\
&+\hbar^3\left(\frac{1}{48} w^3_{xxx}w^{3,4}+\frac{101}{6480} w^3_{xx} w^{3,5}
+\frac{1}{135} w^3_x w^{3,6}+\frac{1}{540} w^3 w^{3,7} \right)\\
&\quad +\frac{\hbar^4}{3600} w^{3,9};\\
\frac{\p w^2}{\p t^{3,0}}=&\frac{1}{36} (w^3)^2 w^2_x-\frac{1}{6} w^2 w^2_x
+\frac{1}{18} w^2 w^3 w^3_x+\frac{1}{2} w^4 w^4_x\\
&+\hbar\left(\frac{1}{12} w^3 w^2_{xxx}+\frac{1}{6} w^2_{xx} w^3_x+\frac{5}{36} w^2_x w^3_{xx}
+\frac{1}{18}  w^2 w^3_{xxx}\right)+\frac{\hbar^2}{20} w^{2,5};\\
\frac{\p w^3}{\p t^{3,0}}=&w^1_x;\\
\frac{\p w^4}{\p t^{3,0}}=&\frac{1}{6} w^4 w^2_x+\frac{1}{18} w^3 w^4 w^3_x
+\frac{1}{36} (w^3)^2 w^4_x+\frac{1}{6} w^2 w^4_x\\
&+\hbar\left(\frac{1}{18} w^4 w^3_{xxx}+\frac{5}{36} w^3_{xx} w^4_x+\frac{1}{6} w^3_x w^4_{xx}
+\frac{1}{12} w^3 w^4_{xxx} \right)+\frac{\hbar^2}{20} w^{4,5};
\end{align*}

\begin{align*}
\frac{\p w^1}{\p t^{4,0}}=&\frac{1}{2} w^4 w^2_x+\frac{1}{6} w^3 w^4 w^3_x
+\frac{1}{12} (w^3)^2 w^4_x+\frac{1}{2} w^2 w^4_x\\
&+\hbar\left( \frac{1}{4} w^3 w^4_{xxx}+\frac{1}{6} w^4 w^3_{xxx}+\frac{5}{12} w^3_{xx} w^4_x
+ \frac12 w^3_x w^4_{xx}\right)+\frac{3\hbar^2}{20} w^{4,5};\\
\frac{\p w^2}{\p t^{4,0}}=&\frac{1}{2} w^3 w^4_x+\frac{1}{2} w^4 w^3_x+\frac{\hbar}2 w^4_{xxx};\\
\frac{\p w^3}{\p t^{4,0}}=& 3 w^4_x;\\
\frac{\p w^4}{\p t^{4,0}}=&w^1_x+\frac{1}{6} w^3 w^2_x+\frac{1}{36} (w^3)^2 w^3_x
+\frac{1}{6} w^2 w^3_x\\
&+\hbar\left( \frac{1}{6} w^2_{xxx}+\frac{1}{36} w^3_x w^3_{xx}+\frac{1}{36} w^3 w^3_{xxx}\right)
+\frac{\hbar^2}{180}  w^{3,5}.
\end{align*}
The flows of the integrable hierarchy can be represented as Hamiltonian systems
\[\frac{\p w^\al}{\p t^{\beta,q}}=\{w^\al(x), H_{\beta,q}\}_1,\]
where the densities $h_{\beta,q}$ of the Hamiltonians $H_{\beta,q}$ are given by 
\begin{equation}\label{zh-11-21}
h_{\beta,q}=\frac{\Gamma(\frac{2\beta-1}6)}{6\Gamma(q+2+\frac{2\beta-1}6)} \left.\res P^{6q+2\beta+5}\right|_{ \p_x^k w^\al\to \left(\frac{\hbar}{2}\right)^{\frac{k}2} \p_x^k w^\al}
\end{equation}
for $\beta=1,2,3,\ q\ge -1$, and 
\begin{equation}
h_{4,q}=\frac{\Gamma(\frac12)}{2\Gamma(q+2+\frac12)} \left.\res Q^{2q+3}\right|_{ \p_x^k w^\al\to \left(\frac{\hbar}{2}\right)^{\frac{k}2} \p_x^k w^\al}
\end{equation}
for $q\ge -1$.
The Hamiltonian structure 
\[ \{w^\al(x),w^\beta(y)\}_1=\eta^{\al\beta} \delta'(x-y)+{\cal O}(\hbar)\]
is a certain deformation of the first Hamiltonian structure of 
hydrodynamic type defined on the loop space of the Frobenius manifold.
As it is shown in \cite{LWZ}, these densities of the Hamiltonians satisfy the tau-symmetry property, and we can obtain from them the two-point functions 
\begin{equation}\label{zh-11-22}
\Omega_{\al,p;\beta,q}(w,w_x,\dots)=\p_x^{-1} \left(\frac{\p h_{\al,p-1}}{\p t^{\beta,q}}\right),\quad \al, \beta=1,\dots,4,\quad p, q\ge 0. 
\end{equation}

The topological solution of the above integrable hierarchy, selected by the string equation, satisfies the following initial condition
\[ \left.w^\al(t)\right|_{t^{\beta,q}=0,\ q\ge 1}=t^{\al,0}.\]
Hence we have
\[ \left.\p_x^k w^\al(t)\right|_{\textrm{small phase space}}=
\left\{\begin{array}{ll} t^{\al,0},&\quad k=0,\\ \delta^\al_1 \delta^k_1, & \quad k\ge 1.\end{array}\right.\]

Now let us construct the tau function by using the following definition: 
\begin{equation}\label{zh-11-7}
\hbar \frac{\p^2\log\tau(t)}{\p t^{\al,p} \p t^{\beta, q}}=\left.\Omega_{\al,p;\beta,q}(w,w_x,\dots)\right|_{w^\al\to {\textrm{topolocal solution}}}.
\end{equation}
To write down the explicit expression of the free energy ${\cal{F}}(t)=\log\tau(t)$, we expand it 
as a power series in $t^{\al,p}\ (\al=1,\dots,n,\ p\ge 1)$ as follows:
\begin{equation}\label{zh-11-10}
\hbar{\cal F}(t)=A(t_0)+\sum_{k\ge 1}
\frac1{k!} \sum_{p_1,\dots, p_k\ge 1} A_{\al_1, p_1;\dots; \al_k, p_k}(t_0) t^{\al_1,p_1}\dots t^{\al_k, p_k}.
\end{equation}
Here and in what follows we indicate the dependence of a function $f$ on 
$t^{1,0},\dots, t^{4,0}$ by writing $f=f(t_0)$. Define
\begin{align*}
&\Omega_{\al,p;\beta,q}(t_0)=\left.\Omega_{\al,p;\beta,q}(w,w_x,\dots)\right|_{\textrm{small phase space}},\\
&\Omega_{\al,p;\beta,q;\gamma, k}(t_0)
=\left.\frac{\p \Omega_{\al,p;\beta,q}(w,w_x,\dots)}{\p t^{\gamma,k}}\right|_{\textrm{small phase space}}.
\end{align*}
In a similar way we define the functions $\Omega_{\al_1,p_1;\al_2,p_2;\dots,\al_k,p_k}(t_0)$
for $k\ge 4$. Then the coefficients  in the expression of the free energy are given by
\begin{align}
&A(t_0)= \int_0^1\sum_{\al=1}^4 \left. t^{\al,0}  \Omega_{1,0;\al,1}(t_0)\right|_{t^{\al,0}\to s\,t^{\al,0}} ds,\label{zh-11-8}\\
&A_{\al,p}(t_0)=\Omega_{\al,p+1;1,0}(t_0),\label{zh-11-9}\\
&A_{\al_1, p_1;\dots; \al_k, p_k}(t_0)=\Omega_{\al_1,p_1;\al_2,p_2;\dots,\al_k,p_k}(t_0).\label{zh-11-12}
\end{align}
Note that the defining relations \eqref{zh-11-7} determine the free energy uniquely up to 
the addition of some constants to the coefficients $A_{\al,p}(t_0)$, and to the addition of some quadratic 
functions of $t^{1,0},\dots, t^{4,0}$ to the leading term $A(t_0)$. To fix these ambiguities, we need to use the following string equation: 
\[ \sum_{p\ge 1} t^{\al,p}\frac{\p {\cal F}(t)}{\p t^{\al,p-1}}+\frac1{2\hbar} \eta_{\al\beta} t^{\al,0}t^{\beta,0}=\frac{\p {\cal F}(t)}{\p t^{1,0}}\]
to obtain the formulae \eqref{zh-11-8}, \eqref{zh-11-9}. Here the matrix $(\eta_{\al\beta})$ is defined by the inverse of the matrix $(\eta^{\al\beta})$
given in \eqref{zh-11-6}.

\begin{rem}
In the expansion formula \eqref{zh-11-10} for the free energy ${\cal F}(t)$, the leading term 
is given by the restriction of ${\cal F}(t)$ to the small phase space. Formula \eqref{zh-11-8}
implies that it equals the restriction of the genus zero free energy ${\cal F}^0(t)$ to the 
small phase space, i.e. 
\[A(t_0)=\left.F(v)\right|_{v^\al\to t^{\al,0}}.\]
We conjecture that for all the simple singularities this property holds true, i.e. the restriction of the 
genus $g\ge 1$ free energy ${\cal F}^g(t)$ to the small phase space vanishes.
\end{rem}

To finish this subsection, let us give some concrete FJRW invariants by using our method.
Due to the limit of computing ability, we have to restrict ourselves to a simple case.
We put $t^{\al,p}=0$ for $p\ge 2$ in the power series expansion \eqref{zh-11-10} of ${\cal F}(t)$ and consider its truncation 
\begin{equation*}\label{zh-11-11}
\left.\hbar{\cal F}\right|_{\textrm{restricted}}(t)=A(t_0)+\sum_{\al} A_{\al,1}(t_0) t^{\al,1}+
\frac1{2!} \sum_{\al,\beta} A_{\al, 1; \beta, 1}(t_0) t^{\al,1} t^{\beta, 1}.
\end{equation*}
By using the formulae given in \eqref{zh-11-8}--\eqref{zh-11-12} we have
\[ \left.\hbar{\cal F}\right|_{\textrm{restricted}}(t)=\tilde{\cal F}^0(t)+\hbar \tilde{\cal F}^1(t),\]
where 
\begin{align*}
&\tilde{\cal F}^0(t)=
\frac{1}{12} {t^{1,0}}({t^{2,0}})^2+\frac{1}{12} ({t^{1,0}})^2 {t^{3,0}}
-\frac{1}{216} ({t^{2,0}})^3{t^{3,0}}+\frac1{1296} ({t^{2,0}})^2({t^{3,0}})^3\\
&+\frac{1}{24} t^{2,0}{t^{3,0}} ({t^{4,0}})^2 +\frac{1}{432} ({t^{3,0}})^3 ({t^{4,0}})^2+\frac{1}{4} {t^{1,0}} ({t^{4,0}})^2-\frac{1}{288} ({t^{2,0}})^4{t^{2,1}}\\
&-\frac{1}{108} ({t^{2,0}})^3{t^{3,0}} {t^{1,1}} -\frac{1}{216} {t^{1,0}} ({t^{2,0}})^3{t^{3,1}}+\frac{1}{48} ({t^{2,0}})^2  ({t^{4,0}})^2 {t^{2,1}}
+\frac{1}{12} {t^{1,0}}({t^{2,0}})^2 {t^{1,1}} \\
&-\frac{1}{72} {t^{1,0}}({t^{2,0}})^2 {t^{3,0}} {t^{2,1}}
+\frac{1}{12} {t^{2,0}} {t^{3,0}} ({t^{4,0}})^2 {t^{1,1}}
+\frac{1}{24} {t^{1,0}}{t^{2,0}}({t^{4,0}})^2 {t^{3,1}}\\
&+\frac{1}{12} ({t^{1,0}})^2 {t^{2,0}}{t^{2,1}}+\frac{1}{12} {t^{2,0}}({t^{4,0}})^3 {t^{4,1}}
+\frac{1}{12} {t^{1,0}} {t^{2,0}} {t^{3,0}} {t^{4,0}} {t^{4,1}} 
+\frac{1}{96}  ({t^{4,0}})^4 {t^{2,1}}\\
&+\frac{1}{4} {t^{1,0}} ({t^{4,0}})^2 {t^{1,1}}+\frac{1}{24} {t^{1,0}}{t^{3,0}} ({t^{4,0}})^2 {t^{2,1}}+\frac{1}{12} ({t^{1,0}})^2 {t^{3,0}}{t^{1,1}} +\frac{1}{36} ({t^{1,0}})^3 {t^{3,1}}\\
&+\frac{1}{4} ({t^{1,0}})^2 {t^{4,0}} {t^{4,1}}+\frac{1}{36}  ({t^{1,0}})^3 ({t^{2,1}})^2
+\frac{1}{12}  ({t^{1,0}})^3 ({t^{4,1}})^2+\frac{1}{18} ({t^{1,0}})^3{t^{1,1}} {t^{3,1}} \\
&+\frac{1}{6}  ({t^{1,0}})^2 {t^{2,0}}{t^{1,1}}{t^{2,1}} 
+\frac{1}{12} ({t^{1,0}})^2 {t^{3,0}} ({t^{1,1}})^2
+\frac{1}{2} ({t^{1,0}})^2 {t^{4,0}}{t^{1,1}}{t^{4,1}}\\
&+\frac{1}{12} {t^{1,0}}({t^{2,0}})^2 ({t^{1,1}})^2 
+\frac{1}{4}{t^{1,0}}({t^{4,0}})^2 ({t^{1,1}})^2+\textrm{156 terms with order higher than 5}.
\end{align*}

\begin{align*}
&\tilde{\cal F}^1(t)=\frac{1}{6} {t^{1,1}}+\frac{1}{432} ({t^{3,0}})^2 {t^{3,1}}
+\frac1{15552} ({t^{3,0}})^4 ({t^{3,1}})^2+\frac{1}{288}({t^{3,0}})^2 ({t^{2,1}})^2 \\
&+\frac{1}{96} ({t^{3,0}})^2 ({t^{4,1}})^2 +\frac{1}{144}({t^{3,0}})^2 {t^{1,1}} {t^{3,1}}+\frac{1}{216} {t^{1,0}} {t^{3,0}} ({t^{3,1}})^2 
 +\frac{1}{12} ({t^{1,1}})^2\\
&+\frac{1}{24} {t^{3,0}} {t^{4,0}} {t^{3,1}}{t^{4,1}} +\frac{1}{72} {t^{2,0}} {t^{3,0}}{t^{2,1}}{t^{3,1}}-\frac{1}{72} {t^{2,0}} ({t^{2,1}})^2
+\frac{1}{432} ({t^{2,0}})^2 ({t^{3,1}})^2\\
&+\frac{1}{144} ({t^{4,0}})^2 ({t^{3,1}})^2
+\frac{1}{24} {t^{2,0}}({t^{4,1}})^2+\frac{1}{12}  {t^{4,0}} {t^{2,1}}{t^{4,1}}.
\end{align*}
By taking derivatives of these $\tilde{\Fcal}^g$ with respect to the time variables $t^{\al,p}$ and then restricting them to the small phase space,
one can obtain all FJRW invariants $\langle \tau_{\alpha_1, p_1}\dots\tau_{\alpha_k, p_k}\rangle$ satisfying
\[\max\{p_i\}\le 1,\quad \#\{i\mid p_i>0\}\le 2.\] 
Note that under the above restriction, there are no invariants with genus greater than one.

\subsection{The $B_3$ Drinfeld-Sokolov hierarchy}

The reduction from the $D_4$ Drinfeld-Sokolov hierarchy to the $B_3$ case can be obtained by taking the invariant locus of $\sigma_1$. According to \eqref{sigma-1}, we only need to take $u^3=0$,
which is equivalent to $\tilde{u}^4=0$ and $w^4=0$.
The super potential and the potential of the Frobenius manifold associated to the Lie algebra of $B_3$ type are obtained from \eqref{zh-11-16} and \eqref{zh-11-3} respectively 
by the substitution $v^4\to 0$. Explicitly, they have the expressions 
\[\lambda(p)=p^6+ v^3 p^4+\left( v^2+\frac14 (v^3)^2\right) p^2+ v^1+\frac16 v^2 v^3+\frac1{108} (v^3)^3.\]
\[F(v)=\frac1{12}v^1(v^2)^2+\frac1{12}(v^1)^2v^3
-\frac{(v^2)^3v^3}{216}+\frac{(v^2)^2(v^3)^3}{1296}
+\frac{(v^3)^7}{1632960}. \]
The associated Drinfeld-Sokolov hierarchy has the following representation:
\begin{equation}\label{zh-11-20}
\frac{\p L}{\p t^{\al,p}}=\frac{\Gamma(\frac{2\al-1}6)}{6 \Gamma(p+1+\frac{2\al-1}6)}
[\left( P^{6 p+2\al-1}\right)_+, L],\quad \al=1,2,3,\quad p\ge 0.
\end{equation}
Here the pseudo-differential operator has the form
\begin{equation}
L=\p_x^6+\sum_{i=1}^3 \p_x^{-1}(\tilde{u}^i \p_x^{2i-1}+\p_x^{2i-1} \tilde{u}^i)
\end{equation}
with $\tilde{u}^1, \tilde{u}^2, \tilde{u}^3$ defined by \eqref{tildu-1}-\eqref{tildu-3}. The pseudo-differential operator
\[P=\p_x+\sum_{i\ge 1}P_i \p_x^{-i}\]
is defined by the equality $P^6=L$.
Then the $B_3$ Drinfeld-Sokolov hierarchy with normalized bilinear form is obtained by performing the substitution \eqref{zh-11-19} to the flows 
\[ \frac{\p w^\al}{\p t^{\beta,q}},\quad \al,\beta=1,2,3,\ q\ge 0\]
defined by \eqref{zh-11-20}.  The densities $h_{\al,p}$ of the Hamiltonians and the functions $\Omega_{\al,p;\beta,q}$ for $\al, \beta=1,2,3, p, q \ge 0$ are defined as in 
\eqref{zh-11-21} and \eqref{zh-11-22}.  Now the free energy ${\cal F}(t)$ is given by
\eqref{zh-11-10}, \eqref{zh-11-9}, \eqref{zh-11-12} and by
\begin{equation}\label{jw-12-28}
A_0(t_0)=\left.F(v)\right|_{v^\al\to t^{\al,0}}.
\end{equation}

Let us compute some concrete invariants. Since the $B_3$ Drinfeld-Sokolov hierarchy is simpler than the $D_4$ one, we can go further now.
Let us take $t^{\al,p}=0$ for $p\ge 3$ in the power series expansion \eqref{zh-11-10} of ${\cal F}(t)$ and consider its truncation 
\begin{align*}
\left.\hbar{\cal F}\right|_{\textrm{restricted}}(t)=&A(t_0)+\sum_{1\le \al\le 3}  \sum_{1\le p\le 2}A_{\al,p}(t_0) t^{\al,p}\\
& + \frac1{2!} \sum_{1\le \alpha, \beta\le 3} \sum_{1\le p, q\le 2} A_{\al, 1; \beta, 1}(t_0) t^{\al,p} t^{\beta, q}.
\end{align*}
By using the formulae given in \eqref{zh-11-9}, \eqref{zh-11-12} and \eqref{jw-12-28},  we have
\[ \left.\hbar{\cal F}\right|_{\textrm{restricted}}(t)=\tilde{\cal F}^0(t)+\hbar \tilde{\cal F}^1(t)+\hbar^2 \tilde{\cal F}^2(t),\]
where $\tilde{\cal F}^0(t)$ has 363 monomial terms which we will not write down here, and 
\begin{align*}
&\tilde{\cal F}^1(t)=\frac1{38880}({t^{3,0}})^5 {t^{3,2}} 
+\frac{1}{432}({t^{3,0}})^3 {t^{1,2}}+\frac{1}{144} {t^{2,0}}({t^{3,0}})^2 {t^{2,2}} +\frac{1}{432} ({t^{3,0}})^2{t^{3,1}} \\
&\quad\quad \quad  +\frac{1}{432} {t^{1,0}}({t^{3,0}})^2{t^{3,2}} 
+\frac{1}{216} ({t^{2,0}})^2 {t^{3,0}}{t^{3,2}} 
+\frac{1}{6} {t^{1,1}}+\frac{1}{6} {t^{1,0}} {t^{1,2}}
-\frac{1}{72} ({t^{2,0}})^2 {t^{2,2}}\\
&\quad\quad \quad+\textrm{85 terms which are quadratic in $t^{\al,1}$ and $t^{\al,2}$}.
\end{align*}
\begin{align*}
&\tilde{\cal F}^2(t)=\frac1{1944}({t^{3,0}})^3 ({t^{3,2}})^2 
+\frac{11}{720}{t^{3,0}} ({t^{2,2}})^2 +\frac{29}{1080}{t^{3,0}}{t^{1,2}} {t^{3,2}} 
+\frac{7}{1080}{t^{1,0}} ({t^{3,2}})^2\\
&\quad\quad \quad+\frac{1}{45}  {t^{2,0}} {t^{2,2}} {t^{3,2}}
+\frac{7}{1080} {t^{3,1}} {t^{3,2}}.
\end{align*}
From these functions, one can obtain all FJRW invariants $\langle \tau_{\alpha_1, p_1}\dots\tau_{\alpha_k, p_k}\rangle$ satisfying
\[\max\{p_i\}\le 2,\quad \#\{i\mid p_i>0\}\le 2.\] 
Note that under the above restriction, there are no invariants with genus greater than two.

\subsection{The $G_2$ Drinfeld-Sokolov hierarchy}

To reduce the $D_4$ Drinfeld-Sokolov hierarchy to the $G_2$ case, we need to choose the lowest weight vectors $\gamma'_1,\dots, \gamma'_4$ given in \eqref{lf-1}, which are eigenvectors of the automorphism $\sigma_4$.
The resulting integrable hierarchy is equivalent to the one defined in Section \ref{sec-51} up to linear coordinate transformations on $u^2$ and $u^3$.
It is easy to see that the invariant locus is obtained by taking $u^2=u^3=0$, which is equivalent to take $w^2=w^4=0$.
Thus the super potential and the potential of the Frobenius manifold associated to the Lie algebra of $G_2$ type are obtained from \eqref{zh-11-16} and \eqref{zh-11-3} respectively 
by the substitution $v^2, v^4\to 0$. They have the expressions 
\[\lambda(p)=p^6+ v^3 p^4+\frac14 (v^3)^2 p^2+ v^1+\frac1{108} (v^3)^3.\]
\[F(v)=\frac1{12}(v^1)^2v^3+\frac{(v^3)^7}{1632960}. \]
The associated Drinfeld-Sokolov hierarchy has the following representation:
\begin{equation}\label{zh-whatever}
\frac{\p L}{\p t^{\al,p}}=\frac{\Gamma(\frac{2\al-1}6)}{6 \Gamma(p+1+\frac{2\al-1}6)}
[\left( P^{6 p+2\al-1}\right)_+, L],\quad \al=1,3,\quad p\ge 0.
\end{equation}
Here the pseudo-differential operator has the form
\begin{equation}
L=\p_x^6+\sum_{i=1}^3 \p_x^{-1}(\tilde{u}^i \p_x^{2i-1}+\p_x^{2i-1} \tilde{u}^i)
\end{equation}
with $\tilde{u}^1, \tilde{u}^2, \tilde{u}^3$ defined by \eqref{tildu-1}-\eqref{tildu-3} with $u^2=u^3=0$. The pseudo-differential operator
\[P=\p_x+\sum_{i\ge 1}P_i \p_x^{-i}\]
is defined by the equality $P^6=L$.
Then the $G_2$ Drinfeld-Sokolov hierarchy with normalized bilinear form is obtained by performing the substitution \eqref{zh-11-19} in the flows 
\[ \frac{\p w^\al}{\p t^{\beta,q}},\quad \al,\beta=1,3,\ q\ge 0\]
defined by \eqref{zh-whatever}.  The densities $h_{\al,p}$ of the Hamiltonians and the functions $\Omega_{\al,p;\beta,q}$ for $\al, \beta=1,3, p, q \ge 0$ are defined as in 
\eqref{zh-11-21} and \eqref{zh-11-22}.  Now the free energy ${\cal F}(t)$ is given by
\eqref{zh-11-10}, \eqref{zh-11-9}, \eqref{zh-11-12} and by
\begin{equation}\label{jw-12-29}
A_0(t_0)=\left.F(v)\right|_{v^\al\to t^{\al,0}}.
\end{equation}

Let us put $t^{\al,p}=0$ for $p\ge 4$ in the power series expansion \eqref{zh-11-10} of ${\cal F}(t)$ and consider its truncation 
\begin{align*}
\left.\hbar{\cal F}\right|_{\textrm{restricted}}(t)=&A(t_0)+\sum_{\al=1,3} \sum_{1\le p\le 3}A_{\al,p}(t_0) t^{\al,1}\\
& +
\frac1{2!} \sum_{\alpha, \beta=1,3} \sum_{1\le p, q\le 3} A_{\al, 1; \beta, 1}(t_0) t^{\al,p} t^{\beta, q}.
\end{align*}
Then by using the formulae given in \eqref{zh-11-9}, \eqref{zh-11-12} and 
\eqref{jw-12-29} we have
\[ \left.\hbar{\cal F}\right|_{\textrm{restricted}}(t)=\tilde{\cal F}^0(t)+\hbar \tilde{\cal F}^1(t)+\hbar^2 \tilde{\cal F}^2(t)+\hbar^3 \tilde{\cal F}^3(t),\]
where $\tilde{\cal F}^0(t)$ has 80 monomial terms and 
\begin{align*}
&\tilde{\cal F}^1(t)=\frac{1}{16796160} (t^{3,0})^8 {t^{3,3}}+\frac{1} {46656} (t^{3,0})^6 {t^{1,3}}+\frac{1}{38880} (t^{3,0})^5 {t^{3,2}}\\
&\quad +\frac{1}{38880}{t^{1,0}}  (t^{3,0})^5 {t^{3,3}}+\frac{1}{432}(t^{3,0})^3 {t^{1,2}}+\frac{1}{432} {t^{1,0}} (t^{3,0})^3{t^{1,3}}\\
&\quad +\frac{1}{432} (t^{3,0})^2 {t^{3,1}}+\frac{1}{432} {t^{1,0}}(t^{3,0})^2 {t^{3,2}}+\frac{1}{864}(t^{1,0})^2 (t^{3,0})^2 {t^{3,3}}\\
&\quad +\frac{1}{6} {t^{1,1}}+\frac{1}{6} {t^{1,0}} {t^{1,2}}+\frac{1}{12}(t^{1,0})^2 {t^{1,3}}\\
&\quad+\textrm{69 terms which are quadratic in $t^{\al,1}, t^{\al,2}$ and $t^{\al,3}$.}
\end{align*}
\begin{align*}
&\tilde{\cal F}^2(t)=\frac{109} {1209323520}(t^{3,0})^9 (t^{3,3})^2+\frac{13}{734832}(t^{3,0})^7 {t^{1,3}} t^{3,3}\\
&
\quad +\frac{13}{839808} {t^{1,0}} (t^{3,0})^6 (t^{3,3})^2
+\frac{13}{839808} (t^{3,0})^6 {t^{3,2}} {t^{3,3}} 
+\frac{61}{62208} (t^{3,0})^5 (t^{1,3})^2\\
&\quad +\frac{31}{20736}(t^{3,0})^4 {t^{1,3}}{t^{3,2}} 
+\frac{73}{62208}(t^{3,0})^4 {t^{1,2}} {t^{3,3}} 
+\frac{83}{31104}{t^{1,0}}(t^{3,0})^4 {t^{1,3}} {t^{3,3}}\\
&\quad  +\frac{1}{1944}(t^{3,0})^3 (t^{3,2})^2 
+\frac{1}{1152}(t^{1,0})^2  (t^{3,0})^3 (t^{3,3})^2
+\frac{11}{15552}(t^{3,0})^3 {t^{3,1}}{t^{3,3}} \\
&\quad +\frac{1}{576} {t^{1,0}}(t^{3,0})^3 {t^{3,2}} {t^{3,3}} 
+\frac{25}{432} {t^{1,0}}(t^{3,0})^2 (t^{1,3})^2 
+\frac{25}{432}(t^{3,0})^2 {t^{1,2}} {t^{1,3}} \\
&\quad +\frac{13}{540} {t^{3,0}}{t^{1,3}} {t^{3,1}}
+\frac{29}{1080}{t^{3,0}} {t^{1,2}} {t^{3,2}} 
+\frac{11}{216} {t^{1,0}}{t^{3,0}} {t^{1,3}} {t^{3,2}}\\
&\quad  +\frac{2}{135} {t^{3,0}}{t^{1,1}} {t^{3,3}} 
+\frac{1}{24} {t^{1,0}} {t^{3,0}} t^{1,2} {t^{3,3}} 
+\frac{5}{108} (t^{1,0})^2 {t^{3,0}}{t^{1,3}} {t^{3,3}}\\
&\quad  +\frac{1}{270}{t^{3,0}} {t^{3,3}} 
+\frac{7}{1080} {t^{1,0}} (t^{3,2})^2
+\frac{5}{1296} (t^{1,0})^3 (t^{3,3})^2
+\frac{7}{1080} {t^{3,1}}{t^{3,2}}\\
&\quad +\frac{11}{1080} {t^{1,0}} {t^{3,1}} {t^{3,3}}
+\frac{5}{432} (t^{1,0})^2 {t^{3,2}} {t^{3,3}}.
\end{align*}

\begin{align*}
&\tilde{\cal F}^3(t)=\frac{7}{3888}(t^{3,0})^2 (t^{3,3})^2
+\frac{281}{9072} {t^{1,3}} {t^{3,3}}.
\end{align*}
From these functions, one can obtain all FJRW invariants $\langle \tau_{\alpha_1, p_1}\dots\tau_{\alpha_k, p_k}\rangle$ satisfying
\[\max\{p_i\}\le 3,\quad \#\{i\mid p_i>0\}\le 2.\] 
Note that under the above restriction, there are no invariants with genus greater than three.
\vskip 0.5truecm 

\noindent{\bf{Acknowledgements}} \ The first and third authors would like to thank Boris Dubrovin for his encouragements and helpful discussions. Part of their work was supported by NSFC No. 11071135, No. 11171176, No. 11222108, and No. 11371214, and by the Marie Curie IRSES project RIMMP. The second author would like to thank his collaborator Todor Milanov from whom he learned most of his knowledge of integrable hierarchies. A special thanks goes to Edward Witten for inspiration and guidance on integrable hierarchies mirror symmetry. His work is partly supported by NSF grant DMS-1103368 and NSF FRG grant DMS-1159265.

\appendix

 \setcounter{equation}{0}
 \renewcommand{\theequation}{A-\arabic{equation}} 

\renewcommand{\thesection}{}

\section{Appendix: Dubrovin-Zhang axiomatic theory of integrable hierarchies}

As we explained in the introduction, Dubrovin-Zhang's approach gives an alternative proof of the   the theorem of Fan-Jarvis-Ruan \cite{FJR1, FJR2}  for the ADE integrable hierarchies conjecture. 
Their proof by-passed Givental's higher genus theory and Kac-Wakimoto hierarchies. In a way, it is  much more direct.  Their proof has been written down explicitly in an extended version of \cite{DZ3} for the $A_n$ case which can be generalized directly 
to the other cases. For the reader's convenience, we present their proof here. 

For any semisimple Frobenius manifold, Dubrovin-Zhang developed an axiomatic theory of integrable hierarchies.  Dubrovin-Zhang's theory starts from a semisimple Frobenius manifold,
to which they associate a bihamiltonian integrable hierarchy of hydrodynamic type, the so-called principal hierarchy. Then, they define the all genus integrable hierarchy as a topological deformation of the principal hierarchy. Under certain
axioms, they are able to show that such a topological deformation is uniquely defined. The ADE Drinfeld-Sokolov hierarchies satisfy these axioms and hence coincide with
Dubrovin-Zhang's topological deformation of the principal hierarchies. The connection to FJRW-theory is through the Virasoro constraints, which are determined by the Frobenius manifold structure. 
More precisely, under the semisimple hypothesis, they show that the genus zero free energy $\Fcal^0$ and the Virasoro constraints uniquely determine the higher
genus free energies $\Fcal^g$. Furthermore,  $e^{\Fcal}$ must be a tau function of their deformed hierarchy and hence a tau function of an ADE Drinfeld-Sokolov hierarchy. 
The only input data from the A-model are (1)  the FJRW Frobenius manifold matches that of the mirror B-model; (2) FJRW-theory satisfies Virasoro constraints. Both
are known for the examples of interest in this paper.
 
\renewcommand{\thesection}{A}

\subsection{The principal hierarchy and the genus zero free energy}\label{liu-1}

Let $\Lambda=\{\Lambda_{g,k}\}$ be a cohomological field theory.  Its genus zero primary potential
(or primary free energy)
\[F(v)=\sum_{k \ge 0}\frac{v^{\alpha_1} \dots v^{\alpha_k}}{k!}\int_{\overline{\Mcal}_{0, k}}\Lambda_{0,k}(\phi_{\alpha_1},\dots,\phi_{\alpha_k})\]
defines the potential of a Frobenius manifold. Namely, 
\begin{itemize}
\item The matrix $\eta_{\alpha\beta}=\p_1\p_{\alpha}\p_{\beta}F$ coincides with the pairing $\langle \phi_{\alpha}, \phi_{\beta}\rangle$.
Here $\p_{\alpha}=\frac{\p}{\p v^{\alpha}}$.
\item The $(1,2)$ tensor $c^{\gamma}_{\alpha\beta}(v)=\eta^{\gamma\gamma'}\p_{\alpha}\p_{\beta}\p_{\gamma'}F$ with $(\eta^{\al\beta})=(\eta_{\al\beta})^{-1}$ defines a family of commutative
associatative products in the following way:
\[\phi_{\alpha}\circ_v \phi_{\beta}=c^{\gamma}_{\alpha\beta}(v)\phi_{\gamma}.\]
We also have
\[\phi_1\circ_v\phi_{\alpha}=\phi_{\alpha}\circ_v\phi_1=\phi_{\alpha},\]
\[\langle \phi_{\alpha}\circ_v\phi_{\beta}, \phi_{\gamma}\rangle=\langle \phi_{\alpha}, \phi_{\beta}\circ_v\phi_{\gamma}\rangle=\p_{\alpha}\p_{\beta}\p_{\gamma}F.\]
\end{itemize}

When $\{\Lambda_{g,k}\}$ comes from the FJRW theory of a quasi-homogeneous
non-degenerate polynomial  $W$, there exists a vector field
\[E=\sum_{\alpha}E^{\alpha}(v)\p_{\alpha}=\sum_{\alpha} d_\alpha v^{\alpha}\p_{\alpha}\]
such that $d_1=1, d_{\alpha}=1-\deg(\phi_{\alpha})$ and $E(F)=(3-c_W)F$. The vector field $E$ is called the Euler vector field of $F$, and $c_W$ is called the charge of $F$.
 If $W$ is an ADE singularity, then $0<d_\alpha \le 1$ for all $1\le \al\le n$, and $0<c_W<1$. From now on, we assume that $F$ comes from the FJRW
theory of an $ADE, D^T_n$  singularity. Then $F$ is always a polynomial, and it coincides with the Frobenius manifold potential obtained from the
mirror Saito Frobenius manifold structure. The latter is known to be isomorphic to the Frobenius manifold of the corresponding Coxeter group \cite{S1, Du2}.

Let us proceed to give the definition of the so-called principal hierarchy of the Frobenius manifold potential $F$ associated to a Coxeter group. For the definition of the principal hierarchy of a general Frobenius manifold, see \cite{Du2, DZ3}.  

First we define $\mu_{\alpha}=1-\frac{1}{2}c_W-d_{\alpha}$, which are called the Hodge grading. Then we introduce  a family of polynomials in $v^\alpha$, denoted by
\[\{\theta_{\alpha,p}(v)\mid \alpha=1, \dots, n,\ p=0, 1, 2, \dots\},\]
via the following relations:
\begin{align*}
&\theta_{\alpha,p}(0)=0,\quad  \theta_{\alpha,0}(v)=\eta_{\alpha\beta}v^{\beta},\quad p\ge 0,\\
&\p_{\alpha}\p_{\beta}\theta_{\gamma,p}(v)=c^{\delta}_{\alpha\beta}(v)\p_{\delta}\theta_{\gamma,p-1}(v),\quad  p\ge 1,\\
&E(\p_{\alpha}\theta_{\beta,p})=(\mu_{\alpha}+\mu_{\beta}+p)\p_{\alpha}\theta_{\beta,p},\quad p\ge 0.
\end{align*}
It is easy to see that the polynomials $\theta_{\alpha,p}$ are uniquely determined by the above conditions, and they  satisfy the normalization condition
\begin{equation}\label{liu-6}
 \p_\xi\theta_\alpha(z) \eta^{\xi\zeta}\p_\zeta\theta_\beta(-z)=\eta_{\al\beta},
\end{equation}
where $\theta_\al(z)=\sum_{p\ge 0} \theta_{\al,p} z^p$.
The set $\{\theta_{\alpha,p}\}$ is called the calibration of $F$.

\begin{rem}
For the readers who are familiar with Givental's notation, the matrix $S(z)$ appearing in Givental's quantization
formula is given
by $S(z)=(S^\alpha_\beta(z))$, where
\[S^\alpha_\beta(z)=\eta^{\alpha\gamma}\sum_{p\ge0}\p_{\gamma} \theta_{\beta,p}z^{-p}.\]
\end{rem}
The principal hierarchy associated to $F$ is a hierarchy of partial differential equations of hydrodynamic type
for $v^1,\dots, v^n$:
\begin{equation}\label{principal}
\frac{\p v^{\beta}}{\p t^{\alpha,p}}=\eta^{\beta\gamma}\p_x\left(\frac{\p \theta_{\alpha,p+1}(v)}{\p v^{\gamma}}\right).
\end{equation}
Since the flow $\frac{\p}{\p t^{1,0}}$ equals $\frac{\p}{\p x}$, we will identify $t^{1,0}$ with $x$ in what follows.
This hierarchy has a bihamiltonian structure 
\begin{equation}\label{liu-17}
\frac{\p v^{\beta}}{\p t^{\alpha,p}}=\{v^{\beta}, H_{\alpha,p}\}_1=\frac1{p+\frac12+\mu_{\alpha}}\{v^{\beta}, H_{\alpha,p-1}\}_2,
\end{equation}
where the Hamiltonians $H_{\alpha,p}=\int h_{\alpha,p}dx$ are defined by the densities
\begin{equation}\label{liu-12-60}
h_{\al,p}=\theta_{\al,p+1},
\end{equation}
and the pair of compatible Poisson 
brackets $\{\ ,\ \}_{1,2}$ are of hydrodynamic type and have the form of the leading terms of \eqref{biham}, i.e.
\begin{align}
&\{v^\al(x), v^\beta(y)\}_1=\eta^{\al\beta}\delta'(x-y), \label{liu-16a}\\
&\{v^\al(x), v^\beta(y)\}_2=g^{\al\beta}(v(x))\delta'(x-y)+\Gamma^{\al\beta}_{\gamma}(v(x)) v^{\gamma}_x(x)\delta(x-y).\label{liu-16b}
\end{align}
Here $(g^{\al\beta})$ is the intersection form of the Frobenius manifold, and $\Gamma^{\al\beta}_\gamma$ are the contravariant Christoffel coefficients of the associated metric \cite{Du2}.
 
The existence of a bihamiltonian structure implies the commutativity
\[\frac{\p}{\p t^{\alpha.p}}\left(\frac{\p v^{\gamma}}{\p t^{\beta.q}}\right)=\frac{\p}{\p t^{\beta,q}}\left(\frac{\p v^{\gamma}}{\p t^{\alpha,p}}\right)\]
of the flows of the principal hierarchy. Hence, all these flows are integrable.
To introduce the tau function of a solution of the principal hierarchy, we first define the functions $\Omega_{\al,p;\beta,q}(v)$ on the Frobenius manifold by
the following generating function \cite{Du2}:
\begin{equation}
\sum_{p, q\ge 0} \Omega_{\al,p;\beta,q}(v) z^p w^q=\frac{\p_\xi\theta_\alpha(z) \eta^{\xi\zeta}\p_\zeta\theta_\beta(w)-\eta_{\al\beta}}{z+w}.
\end{equation}
We call these functions the two-point functions of $F$. They have the following symmetry property:
\[\Omega_{\alpha,p;\beta,q}=\Omega_{\beta,q;\alpha,p}.\]
For any solution $v=v(t)$ of the principal hierarchy we have the identities
\[\frac{\p \Omega_{\alpha,p;\beta,q}(v(t))}{\p t^{\gamma,r}}=\frac{\p \Omega_{\gamma,r;\beta,q}(v(t))}{\p t^{\alpha,p}}.\]
It follows that if $v(t)=\{v^{\alpha}(t)\}$ is a solution to the principal hierarchy \eqref{principal}, then there exists a function $\tau^{0}(t)$ such that
\begin{equation}
\Omega_{\alpha,p;\beta,q}(v(t))=\frac{\p^2 \log \tau^0(t)}{\p t^{\alpha,p}\p t^{\beta,q}}.
\end{equation}
The function $\tau^0(t)$ is called the genus zero tau function corresponding to $v(t)$, and $\Fcal^0(t)=\log \tau^0(t)$ is called a genus zero free energy of $v(t)$.

Note that $\tau^0(t)$ and $\Fcal^0(t)$ are not uniquely defined. If $\Fcal^0(t)$ is a free energy of $v(t)$, then $\tilde{\Fcal}^0(t)=\Fcal^0(t)+\sum_{\alpha,p}c_{\alpha,p}t^{\alpha,p}+c$
is also a free energy of the same solution $v(t)$. The construction given in Section 6 of \cite{Du2} (see also \cite{DZ3}) fixes this ambiguity and obtains the free energy $\Fcal_{\Lambda}^0$ of the cohomological field theory from the topological solutions of the principal hierarchy. 

\begin{pro}[\cite{Du2}]\label{genus-zero}\mbox{}

(a) Define recursively the following sequence in the space of formal power series in $\{t^{\alpha,p}\}$:
\begin{align*}
&v^\beta_{[0]}(t)=t^{\beta,0},\\
&v^\beta_{[k+1]}(t)=\eta^{\beta\gamma}\sum_{\alpha,p}t^{\alpha, p}\left.\frac{\p \theta_{\alpha,p}}{\p v^{\gamma}}\right|_{v^\beta=v^\beta_{[k]}(t)},\quad k=1, 2, \dots.
\end{align*}
Then the limit $v^{\beta}(t)=\lim_{k\to\infty}v^{\beta}_{[k]}(t)$ gives the solution to the Euler-Lagrange equation
\begin{equation}\label{euler-lagr}
v^\beta(t)=\eta^{\beta\gamma}\sum_{\alpha,p}t^{\alpha, p}\frac{\p \theta_{\alpha,p}}{\p v^{\gamma}}.
\end{equation}

(b) Let $v(t)$ be the solution to the above Euler-Laguage equation.  Then it is the unique solution to the principal hierarchy \eqref{principal}
satisfying the initial condition
\begin{equation}\label{init-cond}
v^{\alpha}(t)|_{t^{\alpha,p}=0\ (p>0)}=t^{\alpha,0}.
\end{equation}
This solution is called the topological solution of the principal hierarchy.

(c) The function 
\begin{equation}\label{F0}
\Fcal^0_{\Lambda}=\frac12\sum_{\alpha,p;\beta,q}\Omega_{\alpha,p;\beta,q}(v(t))\tilde{t}^{\alpha,p}\tilde{t}^{\beta,q},\quad \tilde{t}^{\al,p}=t^{\al,p}-\delta^\al_1\delta^p_1
\end{equation}
is a free energy of the topological solution of the principal hierarchy. It coincides with the genus zero free energy of the Frobenius manifold induced from the cohomological field theory $\Lambda$.
\end{pro}

It is shown in \cite{DZ3} that for any Frobenius manifold one can obtain a dense subset of analytic monotonic solutions of the principal hierarchy by solving the following system of equations
\begin{equation}\label{liu-7}
\sum_{p\ge 0} \tilde{t}^{\al,p} \frac{\p\theta_{\al,p}}{\p v^\gamma}=0,\quad \gamma=1,\dots,n,
\end{equation}
where $\tilde{t}^{\al,p}=t^{\al,p}-c^{\al,p}$ and $c^{\al,p}$ are some constants
which vanish except for a finitely many of them. These constants are also required to satisfy a certain genericity condition which is omitted here.  
Note that the above topological solution corresponds to the case when $c^{\al,p}=\delta^\al_1\delta^p_1$. For such a solution of the principal hierarchy, we can fix a free energy by using the formula \eqref{F0}. In what follows we only need to consider such a class of solutions of the principal hierarchy, and we fix their free energies $\Fcal^0$ by using the formula given in  \eqref{F0}. 

\subsection{The topological deformation and the higher genus free energies}\label{liu-2}
As we explained in the last subsection, for any Frobenius manifold one can associate a free energy $\Fcal^0$ to each solution $v(t)=\{ v^\al(t)\}$, obtained by solving the system \eqref{liu-7}, of the principal hierarchy. Now assume that the
Frobenius manifold is semisimple.  Dubrovin-Zhang provided an algorithm in \cite{DZ3} to define a certain topological deformation of the principal hierarchy and the higher genus free energies $\Fcal^g(t) (g\ge 1)$ for solutions of the deformed integrable
hierarchy. 

The main tool  in the construction of Dubrovin-Zhang is the Virasoro symmetries of the principal hierarchy \eqref{principal}. The action of  Virasoro symmetries on the tau functions of the principal hierarchy can be represented by
\[\frac{\p\tau^0}{\p s_m}=a_m^{\al,p;\beta,q} \frac1{\tau^0}\frac{\p\tau^0}{\p t^{\al,p}} \frac{\p\tau^0}{\p t^{\beta,q}}+b^{\beta,q}_{m;\al,p} t^{\al,p}\frac{\p\tau^0}{\p t^{\beta,q}}+c_{m;\al,p;\beta,q} t^{\al,p} t^{\beta,q}\tau^0\]
for $m\ge -1$. Here the coefficients $a, b, c$ define an infinite number of linear operators 
$L_m, m\ge -1$ which satisfy the Virasoro commutation relations
\[ [L_i, L_j]=(i-j)  L_{i+j},\quad i, j\ge -1.\]
We treat the Virasoro operators as a structure associated to a Frobenius manifold \cite{DZ2}. 
For our particular class of
Frobenius manifolds associated to Coxeter groups,  the Virasoro operators have the following expressions:
\begin{align}
L_{-1}=&\sum_{p\ge1}t^{\alpha,p}\frac{\p}{\p t^{\alpha,p-1}}+\frac1{2\hbar}\eta_{\alpha\beta}t^{\alpha,0}t^{\beta,0},\label{liu-19a}\\
L_0=&\sum_{p\ge0}\left(p+\frac12+\mu_{\alpha}\right)t^{\alpha,p}\frac{\p}{\p t^{\alpha,p}}+\frac14\sum_{\alpha=1}^n\left(\frac14-\mu_\alpha^2\right),\label{liu-19b}\\
L_m=&\frac{\hbar}2 \sum_{p+q=m-1} (-1)^{q+1}
\prod_{j=0}^m \left(\mu_\al+j-q-\frac12\right)\eta^{\al\beta} \frac{\p^2}{\p t^{\al,p}\p t^{\beta,q}}\notag\\
&+\quad \sum_{p\ge 0} \prod_{j=0}^m\left(\mu_\al+p+\frac12+j\right) t^{\al,p}\frac{\p}{\p t^{\al,p+m}},\quad m\ge 1.\label{liu-19c}
\end{align}
In particular, they depend only on the pairing (the flat metric) and grading (the Euler vector field) of the Frobenius manifold. 

Note that the actions of the Virasoro symmetries on $\tau^0$ are nonlinear. The idea of \cite{DZ3} is to make a ``change of coordinates" 
\begin{equation}\label{liu-9}
\tau^0(t)\mapsto \tau(t)=e^{\hbar^{-1}\Fcal^0(t)+\Delta F(v, v_x,\dots)|_{v^\al=v^\al(t)}}
\end{equation}
with
\[\Delta F(v,v_x,\dots)=\sum_{g\ge 1}\hbar^{g-1} F^g(v, v_x,\dots),\]
so that in terms of $\tau(t)$ the actions of the Virasoro symmetries are represented as 
\begin{equation}\label{liu-8}
\frac{\p\tau}{\p s_m}=L_m\tau,\quad m\ge -1.
\end{equation}
This linearization condition yields a system of linear equations for the gradients of the
functions $F^g$, which is called the {\em{loop equation}} in \cite{DZ3} (see Theorem 3.10.31 there). It determines
the functions $F^g$ recursively and uniquely up to the addition of constants. Here we assume that $F^g$ depends only on finitely many jet variables, which is a weaker version of Givental's tameness condition \cite{Giv2}. If the function $\Fcal^0(t)$ is the free energy of a solution $v(t)$ given by \eqref{liu-7}, then
the tau function $\tau$ given in \eqref{liu-9} satisfies the Virasoro constraints
\[ L_m|_{t^{\al,p}\to {{t}}^{\al,p}-c^{\al,p}}\tau(t)=0,\quad m\ge -1,\]
and the function
\[\Fcal^g(t)=F^g(v,v_x,\dots)|_{v^\al\to v^\al(t)},\quad g\ge 1\]
is called the genus $g$ free energy associated to $v(t)$. In particular, if we start from the the genus zero free energy $\Fcal_\Lambda^0$ of a cohomological field theory $\Lambda$, then Theorem 3.10.31 of \cite{DZ3} shows that the higher genus free energies $\Fcal_{\Lambda}^g$ are determined uniquely by $\Fcal_\Lambda^0$ and 
the Virasoro constraints 
\[ L_m|_{t^{1,1}\to {{t}}^{1,1}-1}\tau(t)=0,\quad m\ge -1.\]

The functions $F^g(v, v_x,\dots)$ for $g\ge 1$ 
enable us to construct the topological deformation of the principal hierarchy \cite{DZ3}.
 We introduce the new dependent variables
\begin{align}
w^\al=&v^\al+\sum_{g\ge 1} \hbar^g A^\al_g(v, v_x,\dots)\notag\\
=&v^\al+\hbar\eta^{\al\gamma}\frac{\p^2\Delta F(v, v_x,\dots)}{\p t^{1,0}\p t^{\gamma,0}}, \quad \al=1,\dots, n.\label{liu-11}
\end{align}
Such a change of dependent variables is called a {\em quasi-Miura transformation} in \cite{DZ3}. It differs from a Miura-type transformation. In fact, it does not depend {\em polynomially} on the jet variables $v^\al_x,\ 1\le \al\le n$.
Note that the transformation 
\eqref{liu-11} is invertible, so we can also represent $v^\al$ in terms of $w^\beta$ and their jet variables. Now let us write down the principal hierarchy \eqref{principal} in terms of the new variables $w^\al$.  We have an integrable hierarchy of the form
\begin{align}\label{liu-12}
\frac{\p w^\al}{\p t^{\beta,q}}=&\eta^{\beta\gamma}\p_x\left(\frac{\p \theta_{\alpha,p+1}(w)}{\p w^{\gamma}}\right)+\sum_{g\ge 1} \hbar^g R^\al_{g;\beta,q}(w,w_x,\dots,w^{(m_g)}),\\
&\al,\beta=1,\dots,n, \ q\ge 0.\notag
\end{align}
The above hierarchy is just the topological deformation of the principal hierarchy constructed in \cite{DZ3}. 
Since \eqref{liu-11} is a quasi-Miura transformation, it is rather nontrivial that \eqref{liu-12} belongs to the class of integrable hierarchies of KdV-type, i.e. the functions $R^\al_{g;\beta,q}$ depend polynomially on $w^\beta_x, w^\beta_{xx},
\dots,\p_x^{2g+1} w^\beta$ as the Drinfeld-Sokolov hierarchies do. A proof of the polynomial dependence of  $R^\al_{g;\beta,q}$ on $w^\beta_x, w^\beta_{xx}, \dots,\p_x^{2g+1} w^\beta$ is given by Buryak, Posthuma and Shadrin in \cite{BPS1, BPS2}. They also proved the polynomial property of the
first Hamiltonian structure of the deformed principal hierarchy, which is obtained from
the Poisson bracket \eqref{liu-16a} and the first equation of \eqref{liu-17} by 
the change of variables \eqref{liu-11}. In the new variable $w^\al$ the Poisson brackets
$\{\,,\,\}_a, a=1,2$ have the forms
\begin{align}
&\{w^\al(x),w^\beta(y)\}_1=\eta^{\al\beta}\delta'(x-y)\notag\\
&\qquad \qquad +\sum_{g\ge 1}\sum_{k=0}^{2g+1} \hbar^g P^{\al\beta}_{1;g,k}(w,w_x,\dots)\delta^{(2g+1-k)}(x-y),\label{liu-18a}\\
&\{w^\al(x),w^\beta(y)\}_2=g^{\al\beta}(w(x))\delta'(x-y)+\Gamma^{\al\beta}_{\gamma}(w(x)) w^{\gamma}_x(x)\delta(x-y)\notag\\
&\qquad \qquad+\sum_{g\ge 1}\sum_{k=0}^{m_g} \hbar^g P^{\al\beta}_{2;g,k}(w,w_x,\dots)\delta^{(m_g-k)}(x-y).\label{liu-18b}
\end{align}
As proved in \cite{BPS1, BPS2}, the functions $P_{1;g,k}^{\al\beta}$ are homogeneous 
polynomials of $w^\gamma_x, w^\gamma_{xx},\dots$ of degree $k$. The Hamiltonian formalism of the topological deformation of the principal hierarchy is given by
\[\frac{\p w^\al}{\p t^{\beta,q}}=\{w^\al(x),H_{\beta,q}\}_1,\quad \al, \beta=1,\dots, n,\ q\ge 0.\] 
Here the densities of the Hamiltonians $H_{\beta, q}$ are taken as
\[\tilde{h}_{\beta,q}= h_{\beta,q}(v)+\hbar \frac{\p^2\Delta F(v, v_x,\dots)}{\p t^{1,0}\p t^{\beta,q}}
=h_{\beta,q}(w)+\sum_{g\ge 1}\hbar^g Q_{\beta, q; g}(w,w_x,\dots),\]
and $Q_{\beta, q; g}$ are homogeneous polynomials of $w^\gamma_x, w^\gamma_{xx},\dots$ of degree $2g$. Moreover, the functions
\[\tilde{\Omega}_{\al,p;\beta,q}=\hbar \frac{\p^2\log\tau}{\p t^{\al,p}\p t^{\beta,q}}
=\Omega_{\al,p;\beta,q}(w)+\sum_{g\ge 1}\hbar^g W_{\al,p;\beta,q; g}(w,w_x,\dots)\]
provide the tau-structure of the topological deformation of the principal hierarchy. Here $W_{\al,p;\beta,q; g}$ are polynomials of $w^\gamma_x, w^\gamma_{xx},\dots$ of degree $2g$.

It is conjectured in \cite{DZ3} that in the second Poisson bracket \eqref{liu-18b}, the coefficients $P^{\al\beta}_{2;g,k}$ also have the polynomial property.
If this conjecture is valid, then the topological deformation of the principal hierarchy 
has a bihamiltonian structure, a property that is possessed by the KdV hierarchy and, 
more generally, by the Drinfeld-Sokolov hierarchy associated to untwisted affine Lie algebras.

For the semisimple Frobenius manifolds associated to 
the Coxeter groups of ADE type or  ADE singularities, we have a full  picture 
of the topological deformations of the principal hierarchies, including their bihamiltonian structures. Namely, these integrable hierarchies coincide with the Drinfeld-Sokolov hierarchies associated to the untwisted affine Lie algebras of ADE type, which we described in Section \ref{liu-4}. To prove this assertion, we need to use the following results:
\begin{enumerate}
\item The dispersionless ADE Drinfeld-Sokolov hierarchies \eqref{semi-classical} coincide with the principal hierarchies of the Frobenius manifolds associated to the Coxeter groups of ADE type or, equivalently, to the ADE singularities, as it is proved in \cite{DLZ2}. 
\item A result of C.-Z. Wu  \cite{Wu2} shows that  for an ADE Drinfeld-Sokolov hierarchy,
the actions of the Virasoro symmetries on the tau functions defined in \eqref{zh-12-1} are linear. Furthermore,  they are given by  the Virasoro operators \eqref{liu-19a}--\eqref{liu-19c}. 
\item In an extended version of \cite{DZ3}, Dubrovin-Zhang proved that for a semisimple Frobenius manifold, any deformation of the principal hierarchy that possesses a bihamiltonian structure and a tau structure is quasi-trivial.
Moreover, the quasi-Miura transformation has the form \eqref{liu-11}. 
\end{enumerate}
As it is shown in Section \ref{liu-4}, an ADE Drinfeld-Sokolov hierarchy has a 
bihamiltonian structure together with a tau structure. Hence, it is obtained from the principal hierarchy of the Frobenius manifold associated to the corresponding Coxeter group of ADE type by a quasi-Miura transformations of the form \eqref{liu-11}. Now the property of linear actions of the Virasoro symmetries on the tau functions requires that the function $\Delta F$ satisfies the loop equation of the corresponding Frobenius manifold. 
Then using the uniqueness of solution of the loop equation we obtain the following theorem:

\begin{thm}\label{liu-20}
For a semisimple Frobenius manifold associated to a Coxeter group of ADE type, or to 
an ADE singularity, the topological deformation \eqref{liu-12} of the principal hierarchy \eqref{principal} coincides with the ADE Drinfeld-Sokolov hierarchy described in Section \ref{liu-4}.  
\end{thm}

The above theorem provides an alternative proof of  the ADE Witten conjecture by directly connecting the FJRW invariants to the Drinfeld-Sokolov hierarchies.

\begin{thm}[ADE Witten Conjecture]
\begin{description}
\item[(A)] The generating function of the FJRW invariants for an ADE, $D^T_n$  singularity with the maximal diagonal symmetry group is the logarithm of a tau function of the mirror Drinfeld-Sokolov
hierarchy. 
\item[(B)] The generating function of the FJRW invariants for $(D_{2n}, \langle J\rangle)$ is the logarithm of a tau function of the $D_{2n}$-Drinfeld-Sokolov
hierarchy. 
\end{description}
In particular, this tau function is uniquely determined by the Drinfeld-Sokolov hierarchy and the string equation $L_{-1}\tau=0$. 
\end{thm}
\begin{prf}
According to \cite{FJR1, FJR2, FJMR}, the generating function of the genus zero FJRW invariants of an ADE singularity coincides with the genus zero free energy 
of the semisimple Frobenius manifold of the mirror singularity. Furthermore,  the all genus generating function satisfies the Virasoro constraints. By
the result of \cite{DZ3} that the genus zero free energy and the Virasoro constraints 
uniquely determine the full genus free energy, the exponential of this function must be a tau function of the topological deformation of the principal hierarchy
associated to the ADE singularity and satisfies the Virasoro constraints. The theorem then follows from Theorem \ref{liu-20}. 
\end{prf}

\vskip 0.8truecm 
\begin{itemize}
\item[] Department of Mathematical Sciences, Tsinghua University, Beijing 100084, P.R. China\\
\vskip -0.6cm
E-mail address: liusq@mail.tsinghua.edu.cn
\item[] Department of Mathematics, University of Michigan, Ann Arbor, MI 48105 U.S.A\\
\vskip -0.6cm
E-mail address: ruan@umich.edu
\item[] Department of Mathematical Sciences, Tsinghua University, Beijing 100084, P.R. China\\
\vskip -0.6cm
E-mail address: youjin@mail.tsinghua.edu.cn
\end{itemize}
\end{document}